%Format: plain

%% arXive:0805.2277

\input amstex
\documentstyle{amsppt}

\input label.def
\input degt.def
%\input debug.def
%\PrintLabels

{\catcode`\@11
\gdef\proclaimfont@{\sl}}

\newtheorem{Remark}\Remark\thm\endAmSdef
\conjecture\thm\endproclaim

\input epsf
\def\picture#1{$\vcenter{\hbox{\epsffile{#1-bb.eps}}}$}
\def\picture#1{\epsffile{#1-bb.eps}}
\def\plot#1{\centerline{\epsfxsize.8\hsize\picture{plot-#1}}}

\def\dash{\item"\hfill--\hfill"}
\def\dashes{\widestnumber\item{--}\roster}
\def\enddashes{\endroster}

\def\ie{\emph{i.e.}}
\def\eg{\emph{e.g.}}
\def\cf{\emph{cf}}
\def\via{\emph{via}}
\def\etc{\emph{etc}}

\let\Ga\alpha
\let\Gb\beta
\let\Gg\gamma
\let\Gd\delta
\let\Gs\sigma
\def\bGa{\bar\Ga}
\def\bGb{\bar\Gb}
\def\bGg{\bar\Gg}
\def\bGd{\bar\Gd}
\def\bu{\bar u}
\def\bv{\bar v}

\let\eps\epsilon

\def\1{^{-1}}
\def\2{^{-2}}

\def\xt{x_t}
\def\yt{y_t}

\def\bpi{\bar\pi_1}

\loadbold
\def\bA{\bold A}
\def\bD{\bold D}
\def\bE{\bold E}

\let\splus\oplus
\let\bigsplus\bigoplus
\let\onto\twoheadrightarrow
\let\into\hookrightarrow

\def\CG#1{\Z_{#1}}      % the cyclic group
\def\BG#1{\Bbb B_{#1}}  % the braid group
\def\SG#1{\Bbb S_{#1}}  % the symmetric group
\def\AG#1{\Bbb A_{#1}}  % the alternating group
\def\DG#1{\Bbb D_{#1}}  % the dihedral group
\def\RBG#1{\bar\BG{#1}} % the reduced braid group

\def\Cp#1{\Bbb P^{#1}}
\def\B{\bar B}
\def\L{\bar L}
\def\U{\bar U}
\def\tX{\tilde X}
\def\CM{\Cal M}

\def\CK{\Cal K}
\def\term#1-{$\DG{#1}$-}

\def\MB{U} % Milnor ball

\def\PSL{\operatorname{\text{\sl PSL}}}
\def\SL{\operatorname{\text{\sl SL}}}
\def\gcd{\operatorname{g.c.d.}}
\def\discr{\operatorname{discr}}
\def\<#1>{\langle#1\rangle}
\def\ls|#1|{\mathopen|#1\mathclose|}
\def\sdot{{\sssize\bullet}}
\def\DC(#1,#2){\operatorname{Dbl}_{#2}#1}

\def\maplebraid[#1]{{\getmap#1,0,}}
\def\getmap#1,{\afterassignment\getmapp\count0=#1}
\def\getmapp{%\showthe\count0%
 \ifnum\count0<0\let\exp\1\edef\next{-\count0 }\else
 \let\exp\relax\edef\next{\count0 }\fi
 \afterassignment\getmappp\count0=\next}
% \afterassignment\getmappp
% \ifnum\count0<0 \count0=-\count0\else\count1=\count0\fi}
\def\getmappp{%\relax\showthe\count0%
 \ifnum\count0=1\Gd\else
 \ifnum\count0=2\Ga\else
 \ifnum\count0=3\Gb\else
 \ifnum\count0=4\Gg\fi\fi\fi\fi\exp
 \ifnum\count0=0\else\expandafter\getmap\fi}

{\let\1\gdef\let\+\relax
\gdef\setcat{\catcode`\a\active \catcode`\d\active \catcode`\e\active
\catcode`\+\active}
\setcat
\1a[#1]{\bA_{#1}} \1d[#1]{\bD_{#1}} \1e[#1]{\bE_{#1}}
}
{\catcode`\+\active\gdef+{\splus}}
\def\getline{\bgroup\setcat\getlineii}
\def\getlineii[#1,{$#1$\egroup&}
\def\getref"#1"]{\ref{#1}}
\def\tabstrut{\vrule height9.5pt depth2.5pt}

\topmatter

\author
Alex Degtyarev
\endauthor

\title
Fundamental groups of symmetric sextics. II
%The fundamental groups of certain sextics of torus type
\endtitle

\address
Department of Mathematics,
Bilkent University,
06800 Ankara, Turkey
\endaddress

\email
degt\@fen.bilkent.edu.tr
\endemail

\abstract
We study the moduli spaces and compute the fundamental groups of
plane sextics of torus type with the set of inner singularities
$2\bold{A}_8$ or $\bold{A}_{17}$.
We also compute the fundamental groups of a number of other
sextics, both of and not of torus type. The groups found are
simplest possible, \emph{i.e.}, $\Bbb{Z}_2*\Bbb{Z}_3$ and
$\Bbb{Z}_6$, respectively.
\endabstract



\endtopmatter

\document

\section{Introduction}

\subsection{Principal results}
This work concludes the series of papers~\cite{degt.8a2},
\cite{symmetric}, \cite{degt.e6}, where we attempt to classify
and to compute the fundamental groups of irreducible plane sextics
of torus type.
Recall that a sextic~$B$ is said to be of \emph{torus type}
if its equation can be represented in the form $p^3+q^2=0$, where
$p$ and~$q$ are certain homogeneous polynomials of degree~$2$
and~$3$, respectively. Alternatively, $B\subset\Cp2$ is of torus
type if and only if it is the ramification locus of a projection
to~$\Cp2$ of a cubic surface in~$\Cp3$. A representation of the
equation in the form $p^3+q^2=0$ (up to the obvious equivalence)
is called a \emph{torus structure} of~$B$. A singular point~$P$
of~$B$ is called \emph{inner} (\emph{outer}) with respect to a
torus structure~$(p,q)$ if $P$ does (respectively, does not)
belong to the intersection of the conic $\{p=0\}$ and the cubic
$\{q=0\}$.
%The sextic~$B$ is called \emph{tame} if all its
%singular points are inner.
Each sextic~$B$ considered in this paper has a
unique torus structure, see~\cite{degt.Oka}; hence,
we can speak about inner and outer singular points of~$B$. To
indicate the difference, we will use the notation
$(\Sigma_{\text{inner}})\splus\Sigma_{\text{outer}}$ in the
listings. (Note that simple singular points of a sextic are
conveniently identified with their resolution lattices in the
homology of the covering $K3$-surface; for this reason, we use the
direct summation symbol $\splus$ in the notation.)

Another special class is formed by the so called
\emph{\term2n-sextics}, \ie, irreducible plane sextics whose
fundamental group factors to the dihedral group~$\DG{2n}$. Due
to~\cite{degt.Oka}, the \term6-sextics are precisely those of
torus type (see also~\cite{Tokunaga.new}), and the other possible
values are $n=5$ or~$7$. All \term10- and \term14-sextics are
classified and most fundamental groups are computed
in~\cite{degt.Oka3} (see also~\cite{Oka.D10}) and~\cite{degt-Oka}.

First, sextics of torus type appeared in
O.~Zariski~\cite{Zariski.group}.
For the modern state of the subject and further references,
see M.~Oka, D.~T.~Pho~\cite{OkaPho.moduli},
\cite{OkaPho}, H.~Tokunaga \cite{Tokunaga},
and A.~Degtyarev~\cite{degt.Oka}. According
to~\cite{Zariski.group}, the fundamental group
$\pi_1(\Cp2\sminus B)$ of any sextic of torus type factors to the
reduced braid group
$\RBG3:=\BG3/(\Gs_1\Gs_2)^3\cong\PSL(2,\Z)\cong\CG2*\CG3$
(which is the group of the `simplest' curve, the six cuspidal
sextic, constructed in~\cite{Zariski.group}). We show that,
in fact, for
most irreducible sextics of torus type, the group equals~$\RBG3$.
(A summary of known cases is found in Section~\ref{s.torus}. At
present, there are $15$ sets of singularities for which the group
is still unknown.) Our principal results in this paper are
Theorems~\ref{th.moduli} and~\ref{th.group} below, classifying and
computing the fundamental group of sextics with the set of
singularities $(2\bA_8)$ and their degenerations.

\theorem\label{th.moduli}
An irreducible plane sextic of torus type with inner singular
points
$2\bA_8$ or~$\bA_{17}$ has one of the following eight sets of
singularities\rom:
$$
\gathered
(\bA_{17})\splus\bA_2,\quad
 (\bA_{17})\splus\bA_1,\quad
 (\bA_{17}),\\
(2\bA_8)\splus\bA_3,\quad
 (2\bA_8)\splus\bA_2,\quad
 (2\bA_8)\splus2\bA_1,\quad
 (2\bA_8)\splus\bA_1,\quad
 (2\bA_8).\quad
\endgathered
$$
The moduli space of irreducible sextics of torus type realizing
each of the sets of singularities above is unirational\rom; in
particular, it is connected.
\endtheorem

To complete Theorem~\ref{th.moduli},
we also consider reducible sextics of torus type
with a type~$\bA_{17}$ singular point; they split into two cubics.

\theorem\label{th.moduli.3}
A reducible plane sextic of torus type with a type~$\bA_{17}$
singular point has one of the following four
sets of singularities\rom:
$$
(\bA_{17})\splus\bA_2,\quad
 (\bA_{17})\splus2\bA_1,\quad
 (\bA_{17})\splus\bA_1,\quad
 (\bA_{17}).
$$
Each of these sets of singularities is realized by a single connected
deformation family of reducible plane sextics of torus type.
\endtheorem

Theorems~\ref{th.moduli} and~\ref{th.moduli.3} are proved in
Sections~\ref{proof.moduli} and~\ref{proof.moduli.3},
respectively. Another family of reducible sextics of torus type,
those splitting into a quartic and a conic, is considered
in~\S\ref{S.reducible}, see Theorems~\ref{splitting.symmetry}
and~\ref{splitting.moduli}.

\theorem\label{th.group}
The fundamental group $\pi_1(\Cp2\sminus B)$ of each plane sextic~$B$ as
in Theorem~\ref{th.moduli} equals $\RBG3$.
\endtheorem

This theorem is proved in Section~\ref{s.torus.pert}.

\subsection{Other results}
In the proof of Theorems~\ref{th.moduli} and~\ref{th.group}, we
use the approach of~\cite{degt.Oka3}, \cite{degt.8a2},
\cite{degt.e6} (see also Oka~\cite{Oka.symmetric}),
representing the sextics in question
as double coverings of a
certain rigid (maximal in the sense of~\cite{degt.kplets}) trigonal
curve~$\B$ in the Hirzebruch surface~$\Sigma_2$ (see
Section~\ref{s.double}).
According
to~\cite{symmetric}, there are four maximal trigonal curves
admitting a torus structure. Two of them are studied
in~\cite{degt.8a2} and~\cite{degt.e6}, one is considered here
(see~$\B_1$ in Section~\ref{s.a8}), and
the fourth one is reducible (see~$\B_2$ in
Section~\ref{s.a5+a2}). We extend to the remaining reducible
curve~$\B_2$ the results of~\cite{symmetric} and show that it
corresponds to sextics of torus type splitting into a quartic and
a conic, see Theorem~\ref{splitting.symmetry} for the precise
statement. As a consequence, we obtain a deformation
classification of such reducible sextics, see
Theorem~\ref{splitting.moduli}, and compute their fundamental
groups, see~\S\ref{S.groups}. However, we do not make any attempt
to simplify the presentations obtained; we merely summarize the
results
in Theorem~\ref{splitting.group}
and Remark~\ref{rem.group}.

The double covering construction involving the reducible
curve~$\B_2$ makes use of two sections: the linear component
of~$\B_2$ and the ramification locus. Interchanging the sections,
we obtain another family of reducible sextics whose groups are
found with almost no extra work, see~\S\ref{S.others}. The
geometry of these curves is briefly discussed in
Section~\ref{s.remarks}; it involves yet another pair of reducible
maximal trigonal curves found in~\cite{symmetric}.

Instead of simplifying the groups of reducible sextics, we perturb
the curves and obtain the groups of irreducible ones. The
perturbations are constructed using Proposition~5.1.1
in~\cite{degt.8a2}, stating that any induced subgraph of the
combined Dynkin graph os a sextic~$B$ can be realized by a
perturbation of~$B$. This procedure gives rise to a few new
sextics of torus type (the items marked as `see~\ref{S.groups}.?'
in Table~\ref{tab.known} in~\S\ref{S.summary}) and a number of
sextics not of torus type (Table~\ref{tab.abelian}
in~\S\ref{S.summary}). Incorporating the results
of~\cite{degt.8a2} and~\cite{degt.e6}, we obtain the fundamental
groups of all but~$15$ irreducible sextics of torus type (most of
the groups are~$\RBG3$, see Section~\ref{s.torus} for details) and
$768$ other sextics not covered by M.~V.~Nori's
theorem~\cite{Nori} (all groups are abelian). \emph{Extremal} (in
the sense of this paper, \ie, not degenerating to a larger set of
singularities \emph{with known fundamental group}) sets of
singularities with known groups are listed in~\S\ref{S.summary}.

Strictly speaking, for most configurations of singularities, the
connectedness of the equisingular moduli space is still unknown.
For this reason, we state most results in the form of existence
only. However, there is a strong arithmetical evidence for the
following conjecture, which would imply the connectedness.

\conjecture\label{conjecture}
The equisingular moduli space of irreducible plane sextics with
any \emph{non-maximal} configuration of simple singularities is
connected.
%Let~$C_1$ and~$C_2$ be two irreducible plane sextics sharing the
%same configuration of singularities, which are all simple. Assume
%that $C_1$ and~$C_2$ degenerate to the same sextic~$C$. Then $C_1$
%and~$C_2$ are in the same equisingular deformation family.
\endconjecture

Here, the \emph{configuration of singularities} is the set of
singularities enriched with certain information on the mutual
position of the singular points, see~\cite{JAG} for the precise
definition. According to~\cite{degt.Oka}, in the case of
irreducible sextics, all extra information needed is whether the
curve is or is not a \term2n-sextic for some~$n$. A configuration
of singularities is \emph{non-maximal} if it extends to a larger
configuration of singularities still realized by plane sextics.

\subsection{Contents of the paper}
In~\S\ref{S.models}, we explain the double covering construction
used in the proofs, introduce the maximal trigonal curves~$\B_1$
and~$\B_2$, and study the sections of~$\Sigma_2$ that are in a
special position with respect to one of these curves.
Theorem~\ref{th.moduli} is proved here.

\S\ref{S.main} deals with the proof of Theorem~\ref{th.group}. We
sketch out Zariski--van Kampen's method~\cite{vanKampen} in the
special case of
the ruling of~$\Sigma_2$ (section~\ref{s.group}), explain how
the braid monodromy is computed (Section~\ref{s.monodromy}), and
compute the groups of the two maximal sextics
(Sections~\ref{s.a17+a2} and~\ref{s.2a8+a3}). Then, we study the
local perturbations of a few simple singularities
(Section~\ref{s.pert}) and global perturbations of sextics of
torus type (Section~\ref{s.torus.pert}), computing the groups of
the other sextics listed in Theorem~\ref{th.group}.

In~\S\ref{S.groups}, we compute the groups of sextics of torus
type splitting into a quartic and a conic (using the same approach
as in~\S\ref{S.main}). In~\S\ref{S.others}, the representations
obtained are modified to produce the groups of a few other
reducible sextics. In all cases, we are only interested in the
curves whose perturbations contain new irreducible sextics.

\S\ref{S.reducible} is a digression: we establish a geometric
correspondence between the curve $\B_2$ and sextics of torus type
splitting into a quartic and a conic
(Theorem~\ref{splitting.symmetry}). As a consequence, we give a
complete classification of such sextics, see
Theorem~\ref{splitting.moduli}. Theorem~\ref{th.moduli.3} is also
proved here.

In~\S\ref{S.summary}, we give a brief summary of the results
of~\cite{degt.8a2}, \cite{degt.e6}, and this paper.
In particular,
we give a list of the $15$ sets of singularities of sextics of
torus type, for which the fundamental group is still unknown, and
discuss the so called classical Zariski pairs
(Section~\ref{s.classical}).

\subsection{Acknowledgements}
An essential part of this work was completed during my
participation in the special semester on Real and Tropical
Algebraic Geometry held at \emph{Centre Interfacultaire
Bernoulli}, \emph{\'Ecole polytechnique f\'ed\'erale de Lausanne}.
I am thankful to the organizers of the semester and to the
administration of \emph{CIB}.

\section{The trigonal models}\label{S.models}

In Section~\ref{s.double}, we explain the double covering
construction used to produce symmetric plane sextics from trigonal
curves in the Hirzebruch surface~$\Sigma_2$.
Then,
we study two particular rigid trigonal curves~$\B_1$ and~$\B_2$, with the
sets of singularities $\bA_8$ and $\bA_5\splus\bA_2\splus\bA_1$,
respectively.
Most calculations below are straightforward; they were
done using {\tt Maple}.

Note that, according to~\cite{symmetric}, there are four maximal
trigonal curves in~$\Sigma_2$ that admit a torus structure. Two of
them are considered in~\cite{degt.8a2} and~\cite{degt.e6}; the
curves~$\B_1$ and~$\B_2$ studied here are the remaining two.

\subsection{The double covering construction}\label{s.double}
Denote by $\Sigma_2\to\Cp1$ the Hirzebruch surface (\ie, geometrically
ruled rational surface) with an exceptional section~$E$ of
self-intersection~$(-2)$. (One can think of~$\Sigma_2$ as the
minimal resolution of singularities of the quadratic cone
in~$\Cp3$.)
When speaking about affine
coordinates~$(x,y)$ in~$\Sigma_2$, we always assume that $E$ is
given by $y=\infty$.

Any section of~$\Sigma_2$ disjoint from~$E$ has the form
$$
y=s(x):=ax^2+bx+c,\qquad a,b,c\in\C.
\eqtag\label{eq.section}
$$
Given such a section~$\L$, the double covering of the cone
$\Sigma_2/E$ ramified at $E/E$ and~$\L$ is the projective
plane~$\Cp2$, and the deck translation of the covering is an
involutive automorphism $c\:\Cp2\to\Cp2$. Conversely, any
involution $c\:\Cp2\to\Cp2$ has a fixed line~$L_c$ and an isolated
fixed point~$O_c$, and the quotient $\Cp2(O_c)/c$ is the
Hirzebruch surface~$\Sigma_2$. (Here, $\Cp2(O_c)$ stands for the
plane $\Cp2$ blown up at~$O_c$.)

For the purpose of this paper, a \emph{trigonal curve} is a curve
$\B\subset\Sigma_2$ disjoint from~$E$ and intersecting each fiber
at three points. (Alternatively, $\B$ is a curve in $\ls|3E+6F|$
not containing~$E$, where $F$ is a fiber.) Any trigonal curve is
given by a polynomial of the form
$$
f(x,y):=y^3+r_2(x)y^2+r_4(x)y+r_6(x),\qquad\deg r_i=i.
\eqtag\label{eq.curve}
$$
A \emph{torus structure} on the trigonal curve given
by~\eqref{eq.curve} is a decomposition of the form
$$
f(x,y)=\bigl(y+q_2(x)\bigr)^3+\bigl(q_1(x)y+q_3(x)\bigr)^2,
\qquad\deg q_i=i.
\eqtag\label{eq.torus}
$$
A trigonal curve admitting a torus structure is said to be of
\emph{torus type}.

Given a trigonal curve~$\B$ and a section~$\L$ not contained
in~$\B$, the pull-back of~$\B$ under the double covering
$\Cp2\to\Sigma_2/E$ ramified at $E/E$ and~$\L$, see above,
is a plane sextic
$B\subset\Cp2$; we denote it by $\DC(\B,\L)$ and call it the
\emph{double} of~$\B$ ramified at~$\L$. In appropriate affine
coordinates~$(x,y)$ in $\Cp2$, the double is given by the equation
$$
f(x, y^2+s(x))=0,
\eqtag\label{eq.sextic}
$$
where $f$ and~$s$ are as in~\eqref{eq.curve}
and~\eqref{eq.section}, respectively. If $\B$ is of torus type, so
is $\DC(\B,\L)$ for any section~$\L$. The relation between the
singularities of $\B+\L$ and those of $\DC(\B,\L)$ is studied
in~\cite{degt.Oka3}.

The sextic $B=\DC(\B,\L)$ has an involutive \emph{symmetry}, \ie,
an automorphism $c\:\Cp2\to\Cp2$ preserving~$B$. Conversely, any
sextic~$B$ with an involutive symmetry~$c$ such that $O_c\notin B$
is the double of a trigonal curve.

\subsection{The curve~$\B_1$ (the set of singularities $\bA_8$)}\label{s.a8}
The trigonal curve~$\B_1\subset\Sigma_2$ with the set of
singularities $\bA_8$ is a maximal trigonal curve in the sense
of~\cite{degt.kplets}; its skeleton is shown in
Figure~\ref{fig.skeleton}, left. Alternatively, $\B_1$ can be
obtained by a birational transformation from a plane quartic with
a type~$\bA_6$ singular point.
The curve is plotted (in black)
in Figures~\ref{fig.a17+a2} and~\ref{fig.2a8+a3} below.

\midinsert
\centerline{\picture{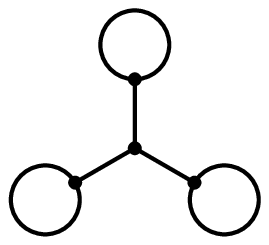}\hfil\picture{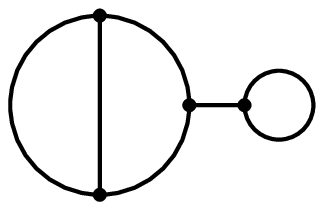}}
\figure
The skeletons of~$\B_1$ and~$\B_2$.
\endfigure\label{fig.skeleton}
\endinsert

In appropriate affine coordinates $(x,y)$ in~$\Sigma_2$, the curve
is given by the polynomial
$$
f_1(x,y):=-y^3+y^2-x^3(2y-x^3).
\eqtag\label{eq.equation}
$$
It has a unique torus structure, given by
$$
f_1(x,y)=(-y)^3+(y-x^3)^2,
$$
and a parametrization
$$
x=x_t:=\frac{t}{t^3+1},\qquad
y=y_t:=\frac1{(t^3+1)^2}.
$$
The discriminant of~$f_1$ with respect to~$y$ is $-x^9(27x^3-4)$.
Thus, $\B_1$ has a type~$\bA_8$ singular point $P_0=(0,0)$
(corresponding to $t=\infty$)
and three vertical tangency points
$$
P_1=\Bigl(\frac{\root3\of4}3,\frac49\Bigr),\quad
P_\pm=\Bigl(\eps_\pm\frac{\root3\of4}3,\frac49\Bigr)
$$
(corresponding to the roots of the equation $2t^3=1$),
where $\eps_\pm=(-1\pm i\sqrt3)/2$ are the primitive cubic roots
of unity.

The surface~$\Sigma_2$ has three automorphisms
$(x,y)\mapsto(\eps x,y)$, $\eps^3=1$,
preserving~$\B_1$ (the \emph{symmetries}
of~$\B_1$)
and three real structures
$\conj_\eps\:(x,y)\mapsto(\eps^2\bar x,\bar y)$, $\eps^3=1$,
with respect to which $\B_1$ is real.
The real part (\ie, the fixed point set)
of $\conj_\eps$ is $(\eps\R,\R)$. In the sequel, we use the
abbreviation
$\conj:=\conj_1$ and $\conj_\pm:=\conj_{\eps_\pm}$.
It is easy to see that a $\conj$-real section~\eqref{eq.section},
$a,b,c\in\R$, $a\ne0$, intersects the real
part $(\eps_\pm\R,\R)$ of~$\conj_\pm$
at two points: $(0,c)$ and
$\Bigl(\eps_\pm\dfrac{b}{a},-\dfrac{b^2}a+c\Bigr)$. (If $a=0$ and
$b\ne0$, the only intersection point is $(0,c)$. If $a=b=0$, the
section is real with respect to all three real structures.)

\subsection{Special sections}\label{s.a8.sections}
Pick a section~$\L$ and consider the double
$B=\DC(\B_1,\L)$. It is a sextic of torus type, see~\ref{s.double}.
If $\L$ is generic, the set of singularities of~$B$ is $(2\bA_8)$;
otherwise, the singularities of~$B$ are recovered from those of
$\B_1+\L$ using the results of~\cite{degt.Oka3}.
Below, for each configuration of $\B_1+\L$,
we parametrize the space of sections~$\L$ that are in the desired
position with respect to~$\B_1$;
for the
reader's convenience, we also indicate the corresponding set
of singularities of~$B$.

A section~\eqref{eq.section} passes through the cusp~$P_0$ of~$\B_1$
(the set of singularities $(\bA_{17})$) if and only if $c=0$.

A section~\eqref{eq.section} is tangent to $\B_1$ at a point
$(\xt(t),\yt(t))$, $2t^3\ne1$,
(the set of singularities $(2\bA_8)\splus\bA_1$)
if and only if $t^3\ne-1$ and
$$
b=-\frac{2t(2t^3-1)a-6t^2}{(2t^3-1)(t^3+1)},\quad
c=\frac{t^2(2t^3-1)a-(4t^3+1)}{(2t^3-1)(t^3+1)^2}
\eqtag\label{eq.tangent}
$$
or $t^3=-1$ and $a=-t$, $b=-2t^2\!/3$. Such a section passes
through the cusp~$P_0$ (the set of singularities
$(\bA_{17})\splus\bA_1$) if and only if
$$
a=\frac{4t^3+1}{t^2(2t^3-1)},\quad
b=-\frac2{t(2t^3-1)},\quad
c=0.
$$

A section~\eqref{eq.section} is inflection tangent to $\B_1$ at a point
$(\xt(t),\yt(t))$, $2t^3\ne1$,
(the set of singularities $(2\bA_8)\splus\bA_2$)
if and only if
$$
a=\frac{3t(8t^6+t^3+2)}{(2t^3-1)^3},\quad
b=-\frac{6t^2(4t^3+1)}{(2t^3-1)^3},\quad
c=\frac{8t^3-1}{(2t^3-1)^3}.
$$
There are three inflection tangents passing through the cusp~$P_0$
of~$\B_1$
(the set of singularities $(\bA_{17})\splus\bA_2$):
$$
t=\frac\eps2,\quad
(a,b,c)=\Bigl(-8\eps,\frac{16\eps^2}3,0\Bigr),
\quad \eps^3=1.
\eqtag\label{eq.inflection}
$$
Clearly, the three tangents~\eqref{eq.inflection}
are interchanged by the
symmetries of~$\B_1$. The tangent~$\L$ corresponding to the real
value
$\eps=1$ is $\conj$-real; it is shown in grey in
Figure~\ref{fig.a17+a2}
below. The tangent intersects~$\B_1$ at~$P_0$
and the following two points:
\dashes
\dash
inflection tangency at $t=\dfrac12$,
$(x,y)=\Bigl(\dfrac{4}{9}, \dfrac{64}{81}\Bigr)\approx(0.44, 0.79)$;
\dash
transversal intersection at $t=-\dfrac32$,
$(x,y)=\Bigl(\dfrac{12}{19}, \dfrac{64}{361}\Bigr)\approx(0.63, 0.78)$.
\enddashes
The intersection of~$\L$ with the real part $\Fix\conj_+$ is at
$(x,y)=\Bigl(-\dfrac{2\eps_+}3,\dfrac{32}9\Bigr)$.

Next lemma deals with the case when a section~$\L$ as
in~\eqref{eq.section} is double tangent
to~$\B_1$ (the set of singularities $(2\bA_8)\splus2\bA_1$).

\lemma\label{double.tangent}
There exists a section~$\L$ of~$\Sigma_2$
tangent to~$\B_1$ at two distinct
points $(\xt(t_i),\yt(t_i))$, $i=1,2$,
$t_1\ne t_2$, $2t_i^3\ne1$, if and only if
$2(t_1+t_2)^3=-1$. A pair $t_1$, $t_2$ as above determines
the double tangent~$\L$ uniquely.
\endlemma

\proof
It suffices to substitute $t=t_1$ and $t=t_2$
to~\eqref{eq.tangent}, equate the resulting expressions for~$b$
and~$c$, and solve for~$a$ the linear system obtained. The
relation $2(t_1+t_2)^3=-1$ is the condition for the compatibility
of the two equations.
\endproof

Letting $t_1=t_2$ in Lemma~\ref{double.tangent}, we obtain three
sections having a point of quadruple intersection with~$\B_1$ (the
set of singularities $(2\bA_8)\splus\bA_3$):
$$
t=\frac\Gd2,\quad
(a,b,c)=\Bigl(-\frac{56\Gd}{27},\frac{64\Gd^2}{81},\frac{256}{243}\Bigr),
\quad\Gd^3=-\frac12.
\eqtag\label{eq.quadruple}
$$
The three sections~\eqref{eq.quadruple}
are interchanged by the symmetries of~$\B_1$.
The section~$\L$ corresponding to the real value $\Gd=-\root3\of4/2$ is
$\conj$-real; it is shown in grey in Figure~\ref{fig.2a8+a3} below.
This section
intersects~$\B_1$ at the following points:
\dashes
\dash
quadruple intersection at $t=-\dfrac{\root3\of4}4$ over
$x=-\dfrac{4\root3\of4}{15}\approx-0.42$;
\dash
transversal intersection at
$t=\Bigl(\dfrac12\pm\dfrac34i\Bigr)\root3\of4$ over
$x=-\Bigl(\dfrac{44}{327}\pm\dfrac{48}{109}i\Bigr)\root3\of4$.
\enddashes
The intersection of~$\L$ with $\Fix\conj_+$ is at
$(x,y)=\Bigl(\dfrac{4\root3\of4\eps_+}{21},\dfrac{512}{567}\Bigr)\approx
 (0.30\eps_+, 0.90)$.
Note that the three points of~$\B_1$ in
%this fiber
the fiber over $x=4\root3\of4\eps_+/21$
are $y\approx0.024$, $0.034$, and $0.94$\strut.

\subsection{Proof of Theorem~\ref{th.moduli}}\label{proof.moduli}
Let $\Sigma$ be a set of singularities as in
Theorem~\ref{th.moduli}, and let $\CM(\Sigma)$ be the moduli space
of irreducible sextics of torus type realizing~$\Sigma$. Due
to~\cite{symmetric}, each sextic~$\B$ in question has a unique
involutive stable symmetry~$c$, and the quotient $B/c$ is a
trigonal curve $\B\subset\Sigma_2$
with a single singular point of type~$\bA_8$.
%isomorphic to~$\B_1$.
Conversely, one has $B=\DC(\B,\L)$ for an appropriate
section~$\L$ (the image of~$L_c$). Hence,
$\CM(\Sigma)$ can be identified with the moduli space of pairs
$(\B,\L)$, where $\B\subset\Sigma_2$ is a trigonal curve with the
set of singularities~$\bA_8$
and $\L$ is a section of~$\Sigma_2$
in a certain prescribed position with respect to~$\B$.
Furthermore, since any curve~$\B$ as above is isomorphic to
the curve~$\B_1$ considered in Section~\ref{s.a8}
(see~\cite{symmetric}) and the group of symmetries of~$\B_1$
is~$\CG3$, there is a cyclic triple ramified covering
$\tilde\CM(\Sigma)\to\CM(\Sigma)$, where $\tilde\CM(\Sigma)$ is
the space of sections~$\L$ forming a prescribed
configuration with~$\B_1$.

The spaces $\tilde\CM(\Sigma)$ are described in
Section~\ref{s.a8.sections}. In each case, $\tilde\CM(\Sigma)$
either is rational
or consists of three rational components (for
$\Sigma=(\bA_{17}\splus\bA_2$, $(2\bA_8)\splus\bA_3$, or
$(2\bA_8)\splus2\bA_1$). In the former case, $\CM(\Sigma)$ is
unirational; in the latter case, the three components are
interchanged by the symmetries of~$\B_1$ (the deck translation of
the covering)
and hence
$\CM(\Sigma)$ is rational.
\qed

\subsection{The curve~$\B_2$
(the set of singularities $(\bA_5\splus\bA_2)\splus\bA_1$)}\label{s.a5+a2}
The trigonal curve $\B_2\subset\Sigma_2$ with the set of
singularities $(\bA_5\splus\bA_2)\splus\bA_1$ is a maximal trigonal curve;
its skeleton is shown in Figure~\ref{fig.skeleton}, right, and the
curve is plotted (in black) in
Figures~\ref{fig.a11+3a2}--\ref{fig.slanted} below.
Alternatively, the curve~$\B_2$ can be obtained by a birational
transformation from a pair of conics inflection tangent to each
other or from a plane quartic splitting into a cuspidal cubic and
a tangent to it.
In appropriate affine coordinates $(x,y)$ in~$\Sigma_2$,
the curve is given by the polynomial
$$
f_2(x,y)=(y^2-x)(y-l(x)),\quad\text{where}\quad
l(x)=-x^2+\frac32x+\frac3{16}.
$$
It splits into a hyperelliptic curve
$\B_2'=\{x=y^2\}$ with a cusp~$R_\infty$ over $x=\infty$ and
a section $\L'=\{y=l(x)\}$. The two components have a point of
inflection tangency at $R_5=(1/4,1/2)$ and a
point of transversal intersection at $R_1=(9/4,-3/2)$,
forming the singular
points of~$\B_2$ of types~$\bA_5$ and~$\bA_1$, respectively.
The only torus structure of~$\B_2$ is given by
$$
64f_2(x,y)=(4y-4x-1)^3+(8xy+6y-12x-1)^2.
$$

Pick a section~$\L$ as in~\eqref{eq.section} and consider the
sextic $B=\DC(\B_2,\L)$. It is of torus type, see~\ref{s.double}.
Furthermore, $B$ splits into a quartic~$B_4$ and conic~$B_2$ (the
pull-backs of~$\B_2'$ and~$\L'$, respectively), which may further
be reducible. If $\L$ is generic, $B_4$ has two cusps and two
points of inflection tangency with~$B_2$.

A section~$\L$  as in~\eqref{eq.section}
is inflection tangent to~$\B_2'$ at a point
$(t^2,t)$, $t\ne0$ (the quartic~$B_4$ has three cusps)
if and only if
$$
(a,b,c)=\Bigl(-\frac1{8t^3},\frac3{4t},\frac{3t}{8}\Bigr).
$$
%as in~\eqref{eq.section}
Such a section passes through~$R_5$ if and only if
$$
%y=27x^2-\frac92x-\frac1{16}.
(a,b,c)=\Bigl(27,-\frac92,-\frac1{16}\Bigr).
\eqtag\label{eq.a11+3a2}
$$
It is shown in grey in Figure~\ref{fig.a11+3a2} below (where $R_1$ is
missing).
The section is inflection tangent to~$\B_2'$ at
%$Q_5=\Bigl(\dfrac1{36},-\dfrac16\Bigr)$
$Q_5=(1/{36},-1/6)$
and intersects~$\L'$ at
%$Q_1=\Bigl(-\dfrac1{28},\dfrac{13}{98}\Bigr)$.
$Q_1=(-1/{28},{13}/{98})$.
The double~$B$ has the set of singularities
$(\bA_{11}\splus2\bA_2)\splus\bA_2\splus2\bA_1$; it splits into a
three cuspidal quartic and a conic.

Making~$\L'$ and~$\L$ trade r\^oles, \ie, considering the sextic
given by
$$
\bigl((y^2+l(x))^2-x\bigr)\bigl(y^2+l(x)-s(x)\bigr)=0,
\eqtag\label{eq.sextic.3}
$$
we obtain a curve $B=\DC((\B_2'+\L),\L')$ with
the set of singularities $(2\bA_5\splus2\bA_2)\splus\bD_5$, also
splitting into a three cuspidal quartic and a conic. (Note that
$B$ is still of torus type, as the trigonal curve
$\B_2'+\L$ is isomorphic to~$\B_2$ \via\ $(x,y)\mapsto(9x,-3y)$.)

A section~$\L$ passes through~$R_5$
and is tangent to~$\B_2'$ at~$R_1$ if and only if
$$
%y=\frac13x^2-\frac{11}6x+\frac{15}{16}.
(a,b,c)=\Bigl(\frac13,-\frac{11}6,\frac{15}{16}\Bigr).
\eqtag\label{eq.a11+2a2+d4}
$$
It is shown in grey in Figure~\ref{fig.a11+2a2+d4} below.
(There is an extra point
%$Q_1=\Bigl(\dfrac{25}4,\dfrac52\Bigr)$
$Q_1=({25}/4,5/2)$
of transversal intersection of~$\L$ and the upper branch of~$\B_2'$;
it is missing in the figure.) The double~$B$ has the set
of singularities $(\bA_{11}\splus2\bA_2)\splus\bD_4$; it splits
into a quartic with the set of singularities $2\bA_2\splus\bA_1$
and a conic.

The section~$\L$ passing through~$R_5$
and tangent to~$\B_2'$ at~$R_\infty$ is given by
$$
y=1/2,%\frac12,
\eqtag\label{eq.a11+e6}
$$
see the solid grey line in Figure~\ref{fig.horizontal}.
The section intersects~$\L'$
transversally at a point $Q_1=(5/4,1/2)$. The double~$B$
has the set of singularities $(\bE_6\splus\bA_{11})\splus2\bA_1$;
it splits into a quartic with a type~$\bE_6$ singular point and a
conic.

The section~$\L$ passing through~$R_1$
and tangent to~$\B$ at~$R_\infty$ is given by
$$
y=-3/2,
\eqtag\label{eq.2a5+e6+a3}
$$
see the dotted grey line in Figure~\ref{fig.horizontal}.
The section intersects~$\L'$
transversally at a point $Q_1=(-3/4,-3/2)$.
The double~$B$
has the set of singularities $(\bE_6\splus2\bA_5)\splus\bA_3$;
it splits into a quartic with a type~$\bE_6$ singular point and a
conic.

A section~$\L$ passes through~$R_\infty$ and is
tangent to~$\B_2'$ at~$R_1$ if and only if
$$
%y=-\frac13x-\frac34,
(a,b,c)=\Bigl(0,-\frac13,-\frac34\Bigr),
\eqtag\label{eq.3a5+d4}
$$
see the solid grey line in Figure~\ref{fig.slanted}.
The section intersects~$\L'$
transversally at the point
%$Q_1=\Bigl(-\dfrac5{12},-\dfrac{11}{18}\Bigr)$.
$Q_1=(-5/{12},-{11}/{18})$.
The double~$B$ has the set of singularities
$(3\bA_5)\splus\bD_4$, splitting into three conics inflection
tangent to each other and having a common point.

A section~$\L$ passes through $P_1$, $P_5$, and~$R_\infty$ if and
only if
$$
%y=-x+\frac34,
(a,b,c)=\Bigl(0,-1,\frac34\Bigr),
\eqtag\label{eq.a11+a5+a3}
$$
see the dotted grey line in Figure~\ref{fig.slanted}.
The double~$B$ has the set of singularities
$(\bA_{11}\splus\bA_5)\splus\bA_3$; it splits into a
quartic with a type~$\bA_5$ singular point and a conic.

\subsection{Other degenerations}\label{s.others}
Here, we consider other possible degeneration of a section~$\L$ with
respect to~$\B_2$, each time showing that the space of sections
admits a rational parametrization. We omit obviously linear
conditions, like passing through one or several of the points~$R_1$,
$R_5$, $R_\infty$.

As above, we fix the notation $B=\DC(\B_2,\L)$ and the splitting
$B=B_4+B_2$ into the pull-backs of~$\B_2'$ and~$\L'$,
respectively.

A section~$\L$ is tangent to~$\B_2'$ at a point $(t^2,t)$,
$t\ne0$, (the quartic~$B_4$ has an extra node) if and
only if
$$
b=-\frac{4at^3-1}{2t},\quad
c=at^4+\frac12t.
$$
Such a section passes through~$R_5$, $R_1$, or~$R_\infty$ if and
only if, respectively,
$$
a=-\frac2{t(2t+1)^2},\quad
a=-\frac2{t(2t-3)^2},\quad\text{or}\quad
a=0.
$$
(The corresponding degenerations of~$B$ are, respectively, the
confluence of two points of inflection tangency of~$B_4$ and~$B_2$
into a single point of $6$-fold intersection, the confluence of
two points of transversal intersection of~$B_4$ and~$B_2$ into a
single tacnode~$\bA_3$, and the confluence of two cusps of~$B_4$
into a single type~$\bA_5$ singular point.)
The section~$\L$ cannot pass through two of the points~$R_5$,
$R_1$, $R_\infty$
unless one of them is the tangency point.

A section~$\L$ is tangent to the $\L'$ component of~$\B_2$
(the conic~$B_2$ splits into two lines) if and only if
$-16ac+3a+4b^2-12b-16c+12=0$. (Clearly, this equation defines a
rational subvariety in the space of triples $(a,b,c)$.) Such a
section cannot pass through~$R_5$ or~$R_1$ (unless it is the
tangency point); it passes through~$R_\infty$ if and only if
$a=0$ and $4c=b^2-3b+3$.

A section tangent to~$B_2'$ at $(t^2,t)$ is also tangent to~$\L'$
if and only if
$$
a=-\frac1{t^2(2t+3)},\quad
b=\frac{3(2t+1)}{2t(2t+3)},\quad
c=\frac{3t}{2(2t+3)}.
$$
When $t\to\infty$, it tends to the section $y=3/4$ tangent
to~$\B_2'$ at~$R_\infty$ and tangent to~$\L'$ (the set of
singularities $(\bE_6\splus2\bA_5)\splus3\bA_1$).

%\section{The fundamental groups: Theorem~\ref{th.group}}

\section{Proof of Theorem~\ref{th.group}}\label{S.main}

In the rest of the paper, we compute the fundamental groups
$\pi_1(\Cp2\sminus B)$ of various sextics~$B$ of the form
$\DC(\B,\L)$, see~\ref{s.double}.
We start with $\B=\B_1$, see~\ref{s.a8}.

Sections~\ref{s.group} and~\ref{s.monodromy} contain an outline of
the approach used to compute the groups. In~\ref{s.a17+a2}
and~\ref{s.2a8+a3}, we study the two maximal doubles of~$\B_1$.
In~\ref{s.pert}, we discuss the perturbations of a few simple
singular points. These results are applied in~\ref{s.torus.pert}
to prove Theorem~\ref{th.group}.

\subsection{Preliminaries}\label{s.group}
Let $\B=\B_1\subset\Sigma_2$ be the trigonal curve as
in~\ref{s.a8}, and let $\L$ be a section of~$\Sigma_2$.
We start with the
group $\bpi:=\pi_1(\Sigma_2\sminus(\B_1\cup\L\cup E))$ and compute it,
applying the classical Zariski--van Kampen method~\cite{vanKampen}
to the ruling of~$\Sigma_2$ (the pencil $\{x=\const\}$ in the notation
of~\S\ref{S.models}).

Pick and fix a real section $S=\{y=\const\gg0\}$ of~$\Sigma_2$ and
a real nonsingular fiber $F=\{x=\tau\}$, where $\tau>0$ is
sufficiently small.
Let $P=F\cap S$, and pick a basis $\Ga$, $\Gb$, $\Gg$, $\Gd$ for the
group $\pi_F:=\pi_1(F\sminus(\B_1\cup\L\cup E),P)$ as shown in
black in
Figure~\ref{fig.basis}, left. (In all cases considered
below, all intersection points are real; the black loops are
oriented counterclockwise.)
Alternatively, denote~$\Ga$, $\Gb$, $\Gg$, $\Gd$ by~$\eta_i$,
$i=1,\ldots,4$, numbering them consecutively according to the
decreasing of the $y$-coordinate. (For example,
in Figure~\ref{fig.basis} one
has $\Ga=\eta_1$, $\Gd=\eta_2$, $\Gb=\eta_3$, $\Gg=\eta_4$.)
We always assume that
$\Ga$, $\Gb$, $\Gg$ are small loops about the three points of
$F\cap\B_1$, numbered consecutively, whereas $\Gd$ is a loop about
$F\cap\L$. Thus, the position of~$\Gd$ in the sequence
$(\Ga,\Gb,\Gg,\Gd)$ may change; this position is important for
some relations.

\midinsert
%\centerline{\picture{vbasis}\hfil\picture{hbasis1}}
\centerline{\valign{\vss\hbox{#}\vss\cr
 \picture{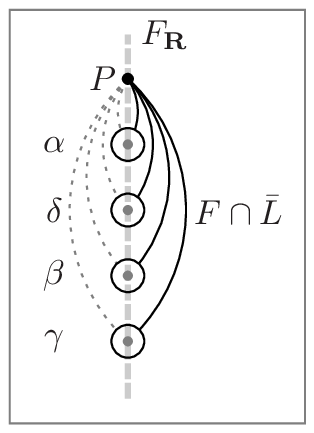}\cr
 \noalign{\hskip3\bigskipamount}\picture{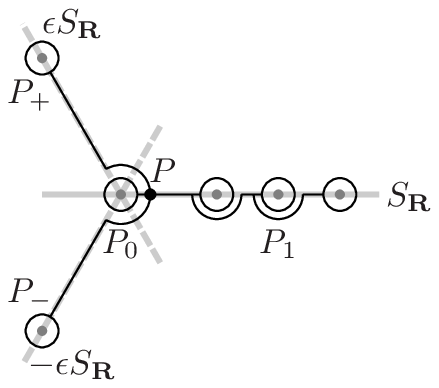}\cr}}
\figure\label{fig.basis}
The basis $\Ga$, $\Gb$, $\Gg$, $\Gd$ and the loops~$\xi_i$
\endfigure
\endinsert

The braid group~$\BG4$ acts on~$\pi_F$: we denote by~$\Gs_1$,
$\Gs_2$, $\Gs_3$ the standard generators of~$\BG4$ and consider
the right action defined by
$\Gs_i\:\eta_i\mapsto\eta_i\eta_{i+1}\eta_i\1$,
$\eta_{i+1}\mapsto\eta_i$.

Let $F_1,\ldots,F_k$ be the singular fibers of $\B_1+\L$.
(Recall that
singular are the fiber $\{x=0\}$ through~$P_0$, the vertical
tangents through~$P_1$ and~$P_\pm$,
and the fibers through the points of
intersection of~$\B_1$ and~$\L$.)
Let $\xi_1,\ldots,\xi_k$ be a basis for the group
$\pi_1(S\sminus(\bigcup F_i\cup\{x=\infty\}),P)$ similar to that
shown in
Figure~\ref{fig.basis}, left: each~$\xi_i$ is a small loop about
$S\cap F_i$ connected to~$P$ by a segment, circumventing the
interfering fibers in the counterclockwise direction. (In the
figures, we consider the section~$\L$ given
by~\eqref{eq.inflection}; necessary modifications for the other
cases are discussed below. The bold grey lines in the figures
represent the real
parts $\Fix\conj_\eps$, $\eps^3=1$.)
For
each~$i$, let
%$m_i\:\pi_F\to\pi_F$
$m_i\in\BG4$
be the \emph{braid monodromy}
along~$\xi_i$, \ie, the automorphism of~$\pi_F$ obtained by
dragging~$F$ along~$\xi_i$ while keeping the base point
on~$\xi_i$. Then, the Zariski--van Kampen theorem~\cite{vanKampen}
states that
$$
\bpi=\bigl<\Ga,\Gb,\Gg,\Gd\bigm|
 \text{$m_i=\id$, $i=1,\ldots,k$,
 $(\eta_1\eta_2\eta_3\eta_4)^2=1$}\bigr>.
\eqtag\label{eq.vanKampen}
$$
Here, each \emph{braid relation} $m_i=\id$ should be understood as
a quadruple of relations $m_i(\Ga)=\Ga$, $m_i(\Gb)=\Gb$,
$m_i(\Gg)=\Gg$, $m_i(\Gd)=\Gd$; the precise form of the
\emph{relation at infinity} $(\eta_1\eta_2\eta_3\eta_4)^2=1$
depends on the
order of the generators.

Now, the passage to the group $\pi_1:=\pi_1(\Cp2\sminus B)$ is
straightforward (see~\cite{degt.Oka3} for details):
one has $\pi_1=\Ker[\kappa\:\bpi/\Gd^2\to\CG2]$, where
$\kappa\:\Ga,\Gb,\Gg\mapsto0$ and $\kappa\:\Gd\mapsto1$. In terms
of the presentations, we
have the following statement.

\lemma\label{bpi->pi}
If $\bpi$ is given by
$\bigl<\Ga,\Gb,\Gg,\Gd\bigm|R_j=1,\ j=1,\ldots,s\bigr>$, then
$$
\pi_1=\bigl<\Ga,\bGa,\Gb,\bGb,\Gg,\bGg\bigm|R'_j=\bar R'_j=1,\
 j=1,\ldots,s\bigr>,
$$
where bar stands for the conjugation by~$\Gd$, $\bar w=\Gd w\Gd$,
and $R'_j$ is obtained from~$R_j$, $j=1,\ldots,s$ by letting
$\Gd^2=1$ and expressing the result in terms of $\Ga$, $\bGa$,
\dots.
\qed
\endlemma

\Remark
Note that $\spbar\:w\mapsto\bar w=\Gd w\Gd$ is an involutive
automorphism of~$\pi_1$. Hence, whenever a relation $R=1$ holds
in~$\pi_1$, the relation $\bar R=1$ also holds.
\endRemark

\subsection{Computing the braid monodromy}\label{s.monodromy}
In this section, we make a few general remarks that facilitate the
computation of the braid monodromy.

According to~\cite{degt.Oka3}, in the presence of the relation at
infinity, (any) one of the braid relations can be ignored. We will
ignore the relation arising from the monodromy around the
cusp~$P_0$.

Since the curves~$\B_1$ and~$\L$, the initial fiber~$F$, and the
base point~$P$ are all chosen $\conj$-real, the
conjugation~$\conj$ induces an automorphism
$\conj_*\:\bpi\to\bpi$. Hence, for each pair $F_\pm$ of conjugate
singular fibers, it suffices to compute the monodromy~$m_+$
about~$F_+$; the relations for~$F_-$ are obtained from those
for~$F_+$ by applying $\conj_*$. The images under $\conj_*$ of the
generators~$\Ga$, $\Gb$, $\Gg$, $\Gd$ are shown in
Figure~\ref{fig.basis}, left: the loops are oriented in the
\emph{clockwise} direction
and connected to~$P$ by the dotted grey paths.
Thus, the action of~$\conj_*$ is as follows:
$$
\gathered
\eta_1\mapsto\eta_1\1,\quad
 \eta_2\mapsto\eta_1\eta_2\1\eta_1\1,\quad
 \eta_3\mapsto(\eta_1\eta_2)\eta_3\1(\eta_1\eta_2)\1,\\
\eta_4\mapsto(\eta_1\eta_2\eta_3)\eta_4\1(\eta_1\eta_2\eta_3)\1.
\endgathered
$$
Its precise form in terms of~$\Ga$, $\Gb$, $\Gg$, $\Gd$ depends on
the order of the generators.

The relations arising from $\conj$-real singular fibers are easily
computed using the plots: the monodromy along a small circle about
the fiber (or along a semicircle circumventing another singular
fiber) is found using a local normal form of the singularity, and
along a segment of the real line the four points of $\B_1+\L$ can
be traced as all but at most two of them are real. In the
computation below, we merely indicate the resulting relations.
(Clearly, it does not really matter whether the interfering
singular fibers are circumvented in the counterclockwise or
clockwise direction; each time, we choose the more convenient
one.)

The monodromy~$m_+$ about~$P_+$ has the form
$\Gs_3^3m_+'\Gs_3^{-3}$, where $m'_+$ is the monodromy along the
small loop about~$P_+$ connected to the point over $x=\eps_+\tau$
by a $\conj_+$-real segment~$I_+$. (For~$m'_+$, we choose the generators
$\Ga'$, $\Gb'$, $\Gg'$ in the fiber~$F'$
over $x=\eps_+\tau$ similar to
Figure~\ref{fig.basis} and take for~$\Gd'$ the image of~$\Gd$
under the monodromy along the arc $x=\tau\exp(it)$,
$t\in[0,2\pi/3]$. Note that the
point $F'\cap\L$ has positive
imaginary part, \ie, in Figure~\ref{fig.basis} it would be located
to the right from $F'_{\R}$.) Now, $m'_+$ is found similar to the
monodromy about~$P_1$, using a plot of the $\conj_+$-real part
of~$\B$, which looks exactly the same as its $\conj$-real part.
However, when computing the braid along~$I_+$, one should take into
account the points of intersection of~$\L$ and the $\conj_+$-real
part $(\eps_+\R,\R)$ of~$\Sigma_2$.

The relations for the conjugate point~$P_-$ are obtained from
those for~$P_+$ by applying $\conj_*$.

\subsection{The set of singularities
$(\bA_{17})\splus\bA_2$}\label{s.a17+a2}
Take for~$\L$ the section given by~\eqref{eq.inflection}.
The curve and the section are plotted in Figure~\ref{fig.a17+a2}.

\midinsert
\plot{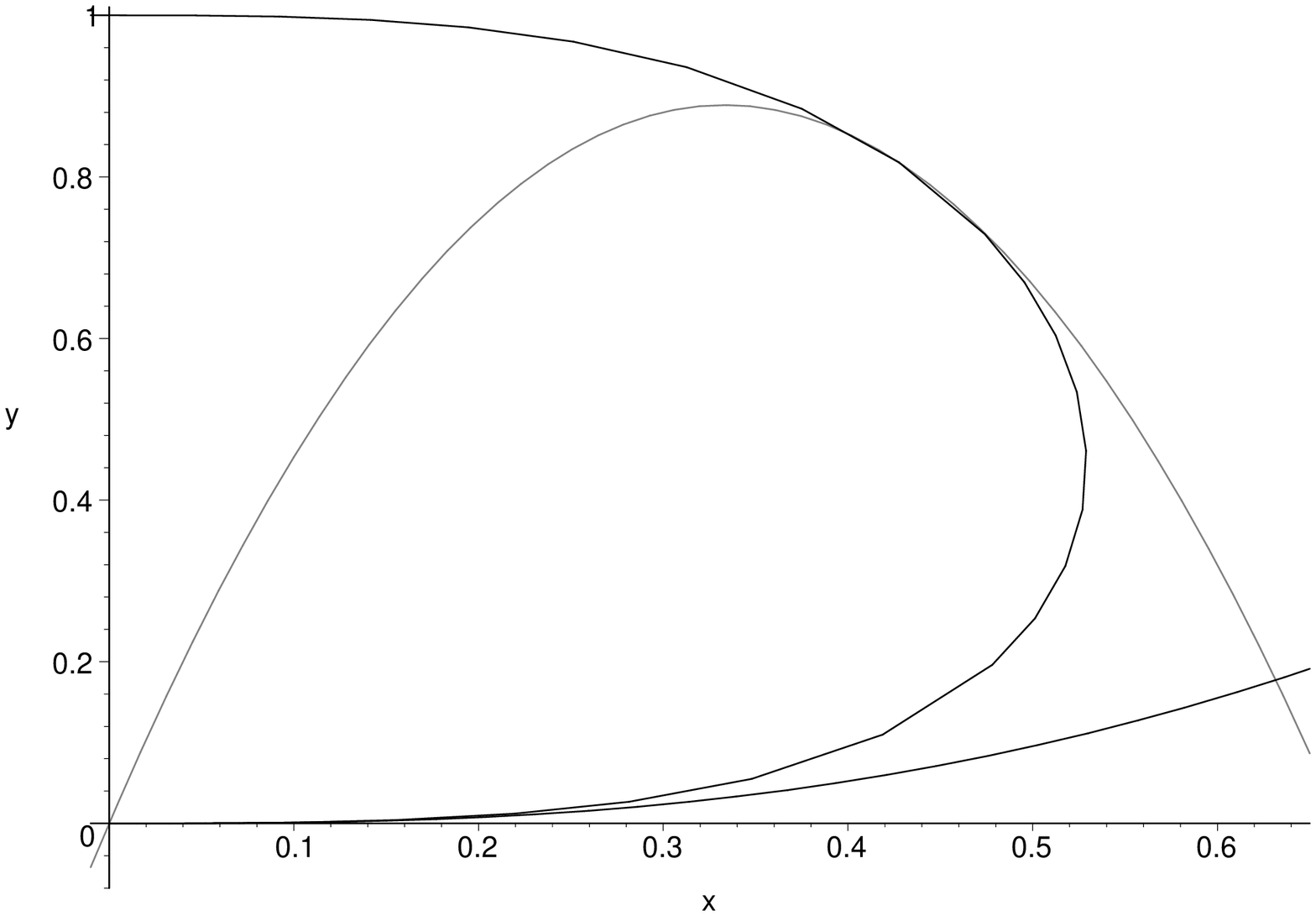}
\figure
The set of singularities $(\bA_{17})\splus\bA_2$
\endfigure\label{fig.a17+a2}
\endinsert

%The complex conjugation induces the automorphism
%$$
%\Ga\mapsto\Ga\1,\quad\Gd\mapsto\Ga\Gd\1\Ga\1,\quad
%\Gb\mapsto(\Ga\Gd)\Gb\1(\Ga\Gd)\1,\quad
%\Gg\mapsto(\Ga\Gd\Gb)\Gg\1(\Ga\Gd\Gb)\1.
%$$
%The relations are:
The basis
%$\Ga$, $\Gd$, $\Gb$, $\Gg$
$(\eta_1,\eta_2,\eta_3,\eta_4)=(\Ga,\Gd,\Gb,\Gg)$
is as shown in Figure~\ref{fig.basis}, and the relations are
$$
\alignat2
&(\Gb\Gg)\1\Gg(\Gb\Gg)=\Ga&\qquad&
 \text{(the vertical tangent through~$P_+$),}\\\allowdisplaybreak
&(\Gd\Gb\Gg\Gb)\Gg(\Gd\Gb\Gg\Gb)\1=\Ga&&
 \text{(the vertical tangent through~$P_-$),}\\\allowdisplaybreak
&(\Ga\Gd)^3=(\Gd\Ga)^3&&
 \text{(the inflection tangency),}\\\allowdisplaybreak
&(\Gd\Ga\Gd)\1\Ga(\Gd\Ga\Gd)=\Gb&&
 \text{(the vertical tangent through~$P_1$),}\\\allowdisplaybreak
&[(\Ga\Gd)\1\Gd(\Ga\Gd),\Gb\Gg\Gb\1]=1&&
 \text{(the transversal intersection),}\\\allowdisplaybreak
%&[\Gd,\Gb\Gg]=1&&
% \text{(???),}\\\allowdisplaybreak
&(\Ga\Gd\Gb\Gg)^2=1&&
 \text{(the relation at infinity).}
\endalignat
$$
(The monodromy~$m_+$ about~$P_+$ is computed as explained in
Section~\ref{s.monodromy}; since $\L$ has no $\conj_+$-real points
over~$I_+$, see Section~\ref{s.a8.sections},
the result is
$\Gs_3^3\Gs_2\Gs_1\Gs_2\1\Gs_3^{-3}$.) Using Lemma~\ref{bpi->pi}, we obtain
the following relations for~$\pi_1$:
$$
\gather
\Gg\Gb\Gg=\Gb\Gg\Ga,\quad
 \bGg\bGb\bGg=\bGb\bGg\bGa,\eqtag\label{eq.1.1}\\\allowdisplaybreak
\Gb\Gg\Gb\Gg=\bGa\Gb\Gg\Gb,\quad
 \bGb\bGg\bGb\bGg=\Ga\bGb\bGg\bGb,\eqtag\label{eq.1.2}\\\allowdisplaybreak
\Ga\bGa\Ga=\bGa\Ga\bGa,\eqtag\label{eq.1.3}\\\allowdisplaybreak
\bGa\1\Ga\bGa=\Gb,\quad
 \Ga\1\bGa\Ga=\bGb,\eqtag\label{eq.1.4}\\\allowdisplaybreak
(\Ga\bGb)\bGg(\Ga\bGb)\1=(\bGa\Gb)\Gg(\bGa\Gb)\1,\eqtag\label{eq.1.5}\\\allowdisplaybreak
%\Gb\Gg=\bGb\bGg,\eqtag\label{eq.1.6}\\\allowdisplaybreak
\Ga\bGb\bGg\bGa\Gb\Gg=1.\eqtag\label{eq.1.6}
\endgather
$$
From~\eqref{eq.1.4} it follows that $\bGa\Gb=\Ga\bGa$ and
$\Ga\bGb=\bGa\Ga$. Substituting these expressions
to~\eqref{eq.1.5}, we obtain
$(\bGa\Ga)\bGg(\bGa\Ga)\1=(\Ga\bGa)\Gg(\Ga\bGa)\1$ or, replacing
$\bGa\Ga$ with $\Ga\bGa\Ga\bGa\1$ from~\eqref{eq.1.3},
$\Ga\1\Gg\Ga=\bGa\1\bGg\bGa$. Thus, introducing
$\Gg_1=\Ga\1\Gg\Ga$ instead of~$\Gg$, we obtain $\bGg_1=\Gg_1$.

Use~\eqref{eq.1.4} and~\eqref{eq.1.3} again to get
$\Ga\1\Gb\Ga=\bGa$ and $\bGa\1\bGb\bGa=\Ga$. Then, the
conjugation of the first and second relations~\eqref{eq.1.1}
by~$\Ga$ and~$\bGa$, respectively, turns them into
$$
\Gg_1\bGa\Gg_1=\bGa\Gg_1\Ga\quad\text{and}\quad
\Gg_1\Ga\Gg_1=\Ga\Gg_1\bGa.
\eqtag\label{eq.1.new}
$$
Similarly, relations~\eqref{eq.1.2}
turn into
$$
\bGa\Gg_1\bGa\Gg_1=\Ga\1(\bGa\Ga\bGa)\Gg_1\bGa
\quad\text{and}\quad
\Ga\Gg_1\Ga\Gg_1=\bGa\1(\Ga\bGa\Ga)\Gg_1\Ga;
$$
using~\eqref{eq.1.3}, they simplify to
$\Gg_1\bGa\Gg_1=\Ga\Gg_1\bGa$ and $\Gg_1\Ga\Gg_1=\bGa\Gg_1\Ga$.
Comparing these expressions to~\eqref{eq.1.new}, one concludes
that $\Ga=\bGa$. Hence, also $\Gb=\bGb=\Ga$ and $\Gg=\bGg$, and
the map $\Ga,\bGa,\Gb,\bGb\mapsto\Gs_1$, $\Gg,\bGg\mapsto\Gs_2$
establishes an isomorphism
$$
\pi_1(\Cp2\sminus B)=\RBG3.
$$
(Relation~\eqref{eq.1.6}, which turns into
$(\Gs_1^2\Gs_2)^2=1$, is equivalent to
$(\Gs_1\Gs_2)^3=1$ in~$\BG3$.)

%\corollary\label{A17->pi}
%The inclusion homomorphism
%$\pi_1(\MB\sminus B)\to\pi_1(\Cp2\sminus B)$ is onto, where $\MB$ is a
%Milnor ball about the type~$\bA_{17}$ singular point of~$B$.
%\qed
%\endcorollary

\subsection{The set of singularities
$(2\bA_8)\splus\bA_3$}\label{s.2a8+a3}
Take for~$\L$ the section given by~\eqref{eq.quadruple}.
The curve and the section are shown in Figure~\ref{fig.2a8+a3}.
\midinsert
\plot{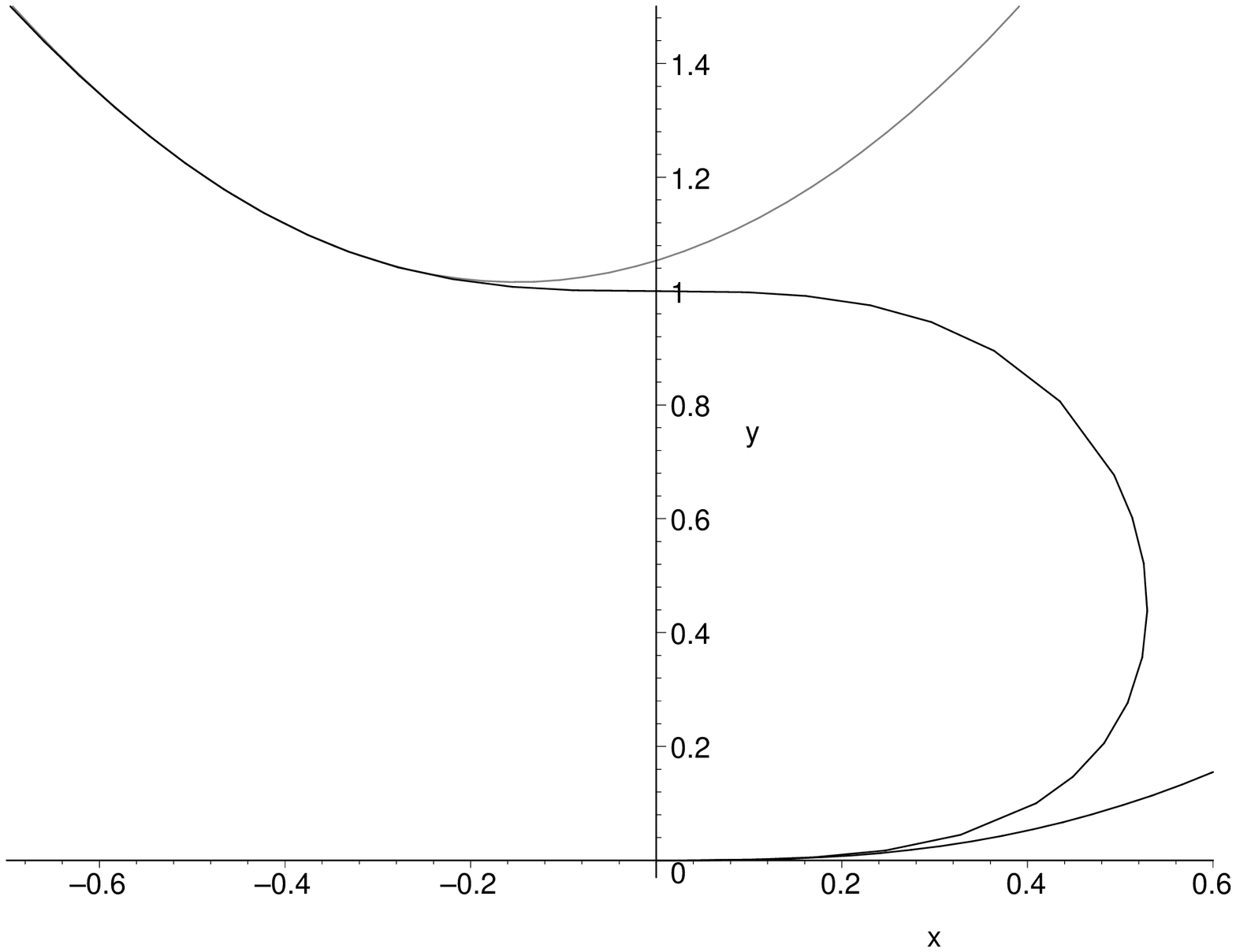}
\figure
The set of singularities $(2\bA_8)\splus\bA_3$
\endfigure\label{fig.2a8+a3}
\endinsert

The singular fibers of $\B_1+\L$ and the basis for
$\pi_1(S\sminus(\bigcup F_i\cup\{x=\infty\}))$ are shown in
Figure~\ref{fig.hbasis}, where $Q_\pm$ are the two points of
transversal intersection of~$\B_1$ and~$\L$, see
Section~\ref{s.a8.sections}.
%The complex conjugation induces the automorphism
%$$
%\Gd\mapsto\Gd\1,\quad\Ga\mapsto\Gd\Ga\1\Gd\1,\quad
%\Gb\mapsto(\Gd\Ga)\Gb\1(\Gd\Ga)\1,\quad
%\Gg\mapsto(\Gd\Ga\Gb)\Gg\1(\Gd\Ga\Gb)\1.
%$$
In the basis $(\eta_1,\eta_2,\eta_3,\eta_4)=(\Gd,\Ga,\Gb,\Gg)$
for the group~$\pi_F$,
the relations are
$$
\alignat2
&(\Gb\Gg)\1\Gg(\Gb\Gg)=(\Gd\Ga)\1\Ga(\Gd\Ga)&\qquad&
 \text{(the vertical tangent through~$P_+$),}\\\allowdisplaybreak
&(\Gb\Gg\Gb)\Gg(\Gb\Gg\Gb)\1=\Gd\Ga\Gd\1&&
 \text{(the vertical tangent through~$P_-$),}\\\allowdisplaybreak
&\Ga=\Gb&&
 \text{(the vertical tangent through~$P_0$),}\\\allowdisplaybreak
&(\Ga\Gd)^4=(\Gd\Ga)^4&&
 \text{(the quadruple intersection point),}\\\allowdisplaybreak
&(\Gd\Ga\Gb\Gg)^2=1&&
 \text{(the relation at infinity).}
\endalignat
$$
The monodromy~$m_+$ about~$P_+$ is found as explained in
Section~\ref{s.monodromy}. This time, $\L$ does have a
$\conj_+$-real point over~$I_+$, see Section~\ref{s.a8.sections},
which is located between the two upper branches of~$\B_1$. Hence,
$m_+=\Gs_3^3\Gs_1^2\Gs_2\Gs_1^{-2}\Gs_3^{-3}$.

\midinsert
\centerline{\picture{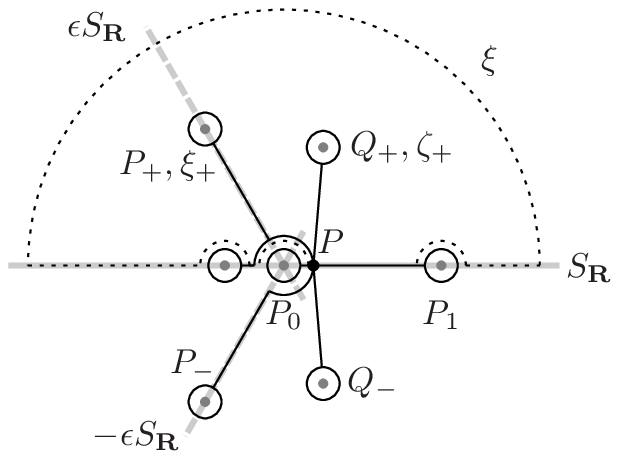}}
\figure
The basis in $\pi_1(S\sminus(\bigcup F_i\cup\{x=\infty\}))$
\endfigure\label{fig.hbasis}
\endinsert

To compute the monodromy~$n_+$ along the loop~$\zeta_+$ about~$Q_+$,
observe that $\zeta_+\xi_+$ is homotopic to the loop~$\xi$
`surrounding' the upper half-plane $\Im x>0$ (the dotted black
loop in Figure~\ref{fig.hbasis}).
More precisely, the loop~$\xi$ is composed of a
large semicircle $x=r\exp(it)$, $t\in[0,\pi]$, $r\gg0$, and
real segment $[-r,r]$ circumventing the real singular fibers in
the \emph{clockwise} direction. The braid over the real segment is
found using local normal forms of the singularities and the plot
of the real part of the curve; it is shown in
Figure~\ref{fig.braid}. The braid over the imaginary
semicircle is the
element $\Delta^2\in\BG4$ corresponding to one full turn. (Indeed,
when $x$ tends to infinity, $\B+\L$ has four quadratic branches
with pairwise distinct leading coefficients.) Hence, the
monodromy~$m$ along~$\xi$ is
$$
m=(\Gs_2\1)(\Gs_3\1\Gs_2)\Delta^2(\Gs_1^{-4})(\Gs_3^{-4}),
$$
and $n_+=mm_+\1$.

\midinsert
\centerline{\picture{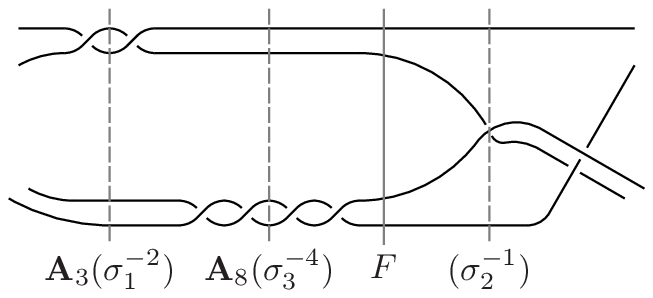}}
\figure
The braid $m_+$ (to be closed by $\Delta^2$)
\endfigure\label{fig.braid}
\endinsert

Equating $n_+(\Gd)=\Gd$, we obtain the relation
$$
%[1, 2, -4, 3, 4, -2, 1, 2, -4, -3, 4, -2, -1]
\Gd\Ga\Gg\1\Gb\Gg\Ga\1\Gd\Ga\Gg\1\Gb\Gg\Ga\1\Gd\1=\Gd,
%\quad\text{or}\quad
%[\Ga\1\Gd\Ga,\Gg\1\Gb\Gg]=1.
$$
resulting from the monodromy about~$Q_+$. The relation from~$Q_-$
is obtained by applying~$\conj_*$. Simplifying both relations, we
get
$$
\gather
[\Ga\1\Gd\Ga,\Gg\1\Gb\Gg]=1,
\eqtag\label{eq.2.+}\\\allowdisplaybreak
[\Gb\1\Gd\Gb,\Gg\Gb\Gg\1]=1.
\eqtag\label{eq.2.-}
\endgather
$$
%$$
%\gather
%\maplebraid[-2, 1, 2, -4, 3, 4, -2, -1, 2, -4, -3, 4, 3, 4, -2, 1, 2, -4, 3, 4, -2, -1, 2, -4, -3, -4, 3, 4, -2, 1, 2, -4, -3, 4, -2, -1, 2]\\
%\maplebraid[-2, 1, 2, -4, 3, 4, -2, -1, 2, -4, -3, 4, 3, 4, -2, 1, 2, -4, -3, 4, -2, -1, 2]
%\endgather
%$$
Eliminate~$\Ga$ using the relation $\Ga=\Gb$,
substitute $\Gb\1\Gd\Gb=\Gd_1$, let $\Gd_1^2=1$, and pass to the
generators $\Gb$, $\Gg$ and $\bGb=\Gd_1\1\Gb\Gd_1$,
$\bGg=\Gd_1\1\Gg\Gd_1$, \cf. Lemma~\ref{bpi->pi}.
We obtain the following set of relations for~$\pi_1$:
$$
\gather
\Gg\Gb\Gg=\Gb\Gg\bGb,\quad
 \bGg\bGb\bGg=\bGb\bGg\Gb,
\eqtag\label{eq.2.1}\\\allowdisplaybreak
\Gg\Gb\Gg=\bGb\Gg\Gb,\quad
 \bGg\bGb\bGg=\Gb\bGg\bGb,
\eqtag\label{eq.2.2}\\\allowdisplaybreak
(\Gb\bGb)^2=(\bGb\Gb)^2,
\eqtag\label{eq.2.3}\\\allowdisplaybreak
%\Gg\1\Gb\Gg=\bGg\1\bGb\bGg,
\Gg\bGb\bGg=\Gb\Gg\Gb,\quad
 \bGg\Gb\Gg=\bGb\bGg\bGb,
\eqtag\label{eq.2.4}\\\allowdisplaybreak
%\Gg\Gb\Gg\1=\bGg\bGb\bGg\1,
\bGg\bGb\Gg=\Gb\Gg\Gb,\quad
 \Gg\Gb\bGg=\bGb\bGg\bGb,
\eqtag\label{eq.2.5}\\\allowdisplaybreak
\Gb\Gg\Gb\bGb\bGg\bGb=1.
\eqtag\label{eq.2.6}
\endgather
$$
Here, relation~\eqref{eq.2.+} turns into
$\Gg\1\Gb\Gg=\bGg\1\bGb\bGg$; in the presence of~\eqref{eq.2.1} it
is equivalent to~\eqref{eq.2.4}. Similarly, \eqref{eq.2.-} turns
into $\Gg\Gb\Gg\1=\bGg\bGb\bGg\1$, which is equivalent
to~\eqref{eq.2.5} in the presence of~\eqref{eq.2.2}.

To simplify the group, consider the generators $u=\Gb\Gg\Gb$ and
$v=\Gg\Gb$, so that $\Gb=uv\1$ and $\Gg=v^2u\1$. Then $\bu=u\1$
(from~\eqref{eq.2.6}) and $\bv=u^2v\2$ (from the first relation
in~\eqref{eq.2.5}). Since the automorphism $w\mapsto\bar w$ is an
involution, one must have $v=\bu^2\bv\2$, \ie,
$$
u^2vu^2=v^2u\2v^2.
\eqtag\label{eq.--=1}
$$
Equate the right hand sides of the first equations
in~\eqref{eq.2.1} and~\eqref{eq.2.2} to obtain
$uvu\2v^2u\2=u\1v^2u\2v$, or $[v,u\2v^2u\2]=1$.
Using~\eqref{eq.--=1}, we get $[v,u^2]=1$; hence $v^3=u^6$ and
$[u,v^3]=1$. Then the first relation in~\eqref{eq.2.1} simplifies
to $u^2=1$; hence also $v^3=1$, $\bu=u$, and $\bv=v$.
%Hence, $\bGb\bGg=\bu\bv\bu\1=u\1v\2u$, and
%substituting this expression to the first relation
%in~\eqref{eq.2.4}, one obtains $u^2=1$. In particular, $\bu=u$ and
%the expression $\bGg=\bv^2\bu\1$ simplifies to $\bGg=v^{-4}u$.
%Substituting the latter to the second relation in~\eqref{eq.2.5},
%one obtains $v^3=1$ and, hence, $\bv=v$.
Thus, $\bGb=\Gb$, $\bGg=\Gg$, and the map $\Gb,\bGb\mapsto\Gs_1$,
$\Gg,\bGg\mapsto\Gs_2$ establishes an isomorphism
$$
\pi_1(\Cp2\sminus B)=\RBG3.
$$

%\corollary\label{A8->pi}
%The inclusion homomorphism
%$\pi_1(\MB\sminus B)\to\pi_1(\Cp2\sminus B)$ is onto, where $\MB$ is a
%Milnor ball about any of the type~$\bA_8$ singular points.
%\qed
%\endcorollary

\subsection{Perturbations of singular points}\label{s.pert}
Consider an isolated singular point~$P$ of a plane curve~$B$, pick
a Milnor ball~$\MB$ about~$P$, and consider a perturbation $B_t$,
$t\in[0,1]$, of $B=B_0$ transversal to $\partial\MB$. Let
$B'=B_1$. We are interested in the perturbation epimorphism
$\pi_1(\MB\sminus B)\onto\pi_1(\MB\sminus B')$,
\cf.~\cite{Zariski.group}.

A perturbation $B$ to~$B'$ is called \emph{maximal} if the total
Milnor number of $B'\cap\MB$ equals $\mu(P)-1$. A perturbation is
called \emph{irreducible} if the normalization of $B'\cap\MB$ is
connected (equivalently, if the abelianization of
$\pi_1(\MB\sminus B')$ is cyclic). If $P$ is simple,
any irreducible perturbation
%of~$P$
factors through a maximal irreducible one.

Let $P$ be of type~$\bA_{3k-1}$ or~$\bE_6$. Then the ball~$\MB$
admits a regular $\SG3$-covering ramified at~$B$. Let
$\phi\:\pi_1(\MB\sminus B)\onto\SG3$ be the corresponding
epimorphism of the fundamental group. A perturbation~$B'$ of~$P$
is said to be of \emph{torus type} if $\phi$ factors through the
perturbation epimorphism above. From the results
of~\cite{degt.Oka} it follows that, if $B$ is a sextic of torus
type and $P$ is an inner singularity of~$B$, then $B'$ is still of
torus type if and only if the perturbation is of torus type.

According to E.~Looijenga~\cite{Looijenga},
the deformation classes (in the obvious sense) of perturbations of
a simple singularity~$P$ are enumerated by the induced subgraphs
of the Dynkin diagram of~$P$ (up to a certain equivalence, which
is not important here). In the statements and proofs below, we
merely indicate the result of the enumeration.

%\lemma\label{pert.a17}
%For any perturbation of a type~$\bA_{17}$ singular point that is
%not of torus type, the group $\pi_1(\MB\sminus B')$ is abelian.
%\endlemma
%
%\lemma\label{pert.a8}
%For any perturbation of a type~$\bA_8$ singular point that is
%not of torus type, the group $\pi_1(\MB\sminus B')$ is cyclic.
%\endlemma

\lemma\label{pert.a11}
The maximal irreducible perturbation of a type~$\bA_{11}$ singular
point are
$\bA_8\splus\bA_2$ \rom(of torus type\rom) and
$\bA_{10}$ and $\bA_6\splus\bA_4$ \rom(not of torus type\rom).
In the former case, the group $\pi_1(\MB\sminus B')$ equals $\BG3$\rom;
in the latter case, it is cyclic.
\endlemma

\lemma\label{pert.a5}
The maximal irreducible perturbation of a type~$\bA_5$ singular
point are
$2\bA_2$ \rom(of torus type\rom) and
$\bA_4$ \rom(not of torus type\rom).
In the former case, the group $\pi_1(\MB\sminus B')$ equals $\BG3$\rom;
in the latter case, it is cyclic.
\endlemma

\lemma\label{pert.a7}
The maximal irreducible perturbation of a type~$\bA_7$ singular
point are $\bA_6$ and $\bA_4\splus\bA_2$.
For these perturbations,
the group $\pi_1(\MB\sminus B')$ is cyclic.
\endlemma

\proof[Proof of Lemmas~\ref{pert.a11}--\ref{pert.a7}]
All statements are a well known property of type~$\bA$ singular
points: any perturbation of a type~$\bA_p$ singular point has the
set of singularities $\bigsplus\bA_{p_i}$ with
$d=(p+1)-\sum(p_i+1)\ge0$, and the group $\pi_1(\MB\sminus B')$ is
given by $\<\Ga,\Gb\,|\,\Gs^s\Ga=\Ga,\ \Gs^s\Gb=\Gb>$, where $\Gs$
is the standard generator of the braid group~$\BG2$ acting on
$\<\Ga,\Gb>$ and $s=1$ if $d>0$ or $s=\gcd(p_i+1)$ if $d=0$.

A perturbation is reducible if and only if $s$ above is even; a
perturbation is of torus type is and only if $s=0\bmod3$.
\endproof

\lemma\label{pert.d5}
The only maximal irreducible perturbation of a type~$\bD_5$
singular point is~$\bA_4$.
For this perturbation,
the group $\pi_1(\MB\sminus B')$ is cyclic.
\endlemma

\proof
The perturbation~$B_t$ can be realized by a family $C_t\subset\C^2$
of affine
quartics with a point of quadruple intersection with the line at
infinity, so that
$(\MB,B_t)\cong(\C^2,C_t)$ for each $t\in[0,1]$. For a
quartic~$C_1$ with a type~$\bA_4$ singular point, one has
$\pi_1(\C^2\sminus C_1)=\Z$, see~\cite{groups}.
\endproof

\lemma\label{pert.d4}
For any nontrivial perturbation of a type~$\bD_4$
singular point, the group $\pi_1(\MB\sminus B')$ is abelian.
\endlemma

\proof
Any perturbation~$B_t$ of a type~$\bD_4$ singular point can be
realized by a family $C_t\subset\C^2$ of affine cubics transversal
to the line at infinity, so that $(\MB,B_t)\cong(\C^2,C_t)$ for
each $t\in[0,1]$. Unless $C_t$ is a triple of lines passing
through a single point, the group $\pi_1(\C^2\sminus C_t)$ is
abelian.
\endproof

\lemma\label{pert.e7}
The maximal irreducible perturbations of a type~$\bE_7$ singular
point are $\bE_6$, $\bA_6$, and $\bA_4\splus\bA_2$.
For the perturbation $\bA_4\splus\bA_2$, one has
$\pi_1(\MB\sminus B')=\Z\times\SL(2,\Bbb F_5)$\rom;
for other irreducible perturbations,
the group $\pi_1(\MB\sminus B')$ is cyclic.
\endlemma

\proof
Any perturbation~$B_t$ of a type~$\bE_7$ singular point can be
realized by a family $C_t\subset\C^2$ of affine quartics
inflection tangent to the line at infinity (see,
\eg,~\cite{quintics}),
so that $(\MB,B_t)\cong(\C^2,C_t)$ for
each $t\in[0,1]$.
The groups $\pi_1(\C^2\sminus C_t)$ for such
quartics are found in~\cite{groups}.
\endproof

We need an explicit description of the epimorphism
$\pi_1(\MB\sminus B)\onto\pi_1(\MB\sminus B')$ for the
perturbation $\bE_7\to\bA_4\splus\bA_2$. Consider the isotrivial
trigonal curve $\B\subset\Sigma_2$ given by $y^3+x^3y=0$. It has
two singular fibers: a fiber of type~$\tilde\bE_7$ at $P=(0,0)$ and
a fiber of type~$\tilde\bA_1^*$ over $x=\infty$. Let $\U$ be the
affine part of~$\Sigma_2$ (the complement of the exceptional
section~$E$ and the fiber over $x=\infty$). Then the pair
$(\U,\B)$ is diffeomorphic to the pair $(U,B)$ above, and the
perturbation can be realized by deforming~$\B$ to a maximal
trigonal curve~$\B'$ with singular fibers of types~$\tilde\bA_4$,
$\tilde\bA_2$, $\tilde\bA_0^*$, and~$\tilde\bA_1^*$ (keeping the
last one at infinity).

\midinsert
\centerline{\picture{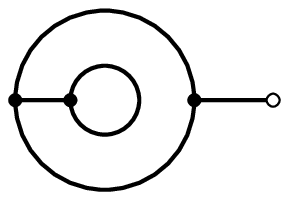}}
\figure
The skeleton of~$\B'$
\endfigure\label{fig.a4+a2}
\endinsert

The skeleton of the perturbed curve~$\B'$ is shown in
Figure~\ref{fig.a4+a2}. It follows that $\B'$ and all its real
fibers can be chosen real, and one can use the skeleton to sketch
the real part of~$\B'$ and compute the braid monodromy.
(Alternatively, one can use the description of the braid monodromy
in terms of the skeleton given in~\cite{degt.kplets}.) As a
result, in the standard generators $\Ga$, $\Gb$, $\Gg$, \cf.
Figure~\ref{fig.basis}, in the fiber over $x\ll0$ (in which
both~$\B$ and~$\B'$ have three real points), the relations for
$\pi_1(\U\sminus\B')$ are
$$
\Ga\Gb\Ga=\Gb\Ga\Gb,\quad
\Ga\Gg\Ga=\Gg\Ga\Gb,\quad
\Gg\Ga\Gg=\Gb\Gg\Ga.
\eqtag\label{eq.e7}
$$

\subsection{Perturbations of sextics of torus type}\label{s.torus.pert}
Here, we state a few simple lemmas about irreducible plane
sextics of
torus type with the fundamental group
$\pi_1(\Cp2\sminus B)=\RBG3$. (Note that, in fact,
any sextic whose group is $\RBG3$ is irreducible and
then, due to~\cite{degt.Oka}, it is of torus type.)

%\lemma\label{RBG->torus}
%Let $B$ be a plane sextic of torus type. Then any
%epimorphism $\phi\:\RBG3\onto\pi_1(\Cp2\sminus B)$ is
%an isomorphism.
%\endlemma
%
%\proof
%There are epimorphisms
%$\RBG3\onto\pi_1(\Cp2\sminus B')\onto\RBG3$, where the first arrow
%is~$\phi$ and the last one is induced by the perturbation of~$B$
%to Zariski's six cuspidal sextic, see~\cite{Zariski.group}. Since
%the group $\RBG3=\PSL(2,\Z)$ is obviously residually finite, hence
%Hopfian, both epimorphisms are isomorphisms.
%\endproof

\lemma\label{RBG->torus}
Let $B$ be an irreducible plane sextic of torus type. Then any
epimorphism $\phi\:\BG3\onto\pi_1(\Cp2\sminus B)$ factors through
an isomorphism $\RBG3=\pi_1(\Cp2\sminus B)$.
\endlemma

\proof
Recall that $(\Gs_1\Gs_2)^3\subset\BG3$ is a central element whose
image in the abelianization $\BG3/[\BG3,\BG3]=\Z$ is~$6$. Since
the abelianization of $\pi_1(\Cp2\sminus B)$ is $\CG6$, one has
$(\Gs_1\Gs_2)^3\in\Ker\phi$ and $\phi$ factors to an epimorphism
$\RBG3\onto\pi_1(\Cp2\sminus B)$.
On the other hand, there is an epimorphisms
$\pi_1(\Cp2\sminus B')\onto\RBG3$ induced by the perturbation of~$B$
to Zariski's six cuspidal sextic, see~\cite{Zariski.group}. Since
the group $\RBG3=\PSL(2,\Z)$ is obviously residually finite, hence
Hopfian, both epimorphisms are isomorphisms.
\endproof

\corollary\label{torus->torus}
Let~$B$ be a sextic of torus type,
$\pi_1(\Cp2\sminus B)=\RBG3$,
and let~$B'$ be a perturbation of~$B$ which is also of torus type.
Then $\pi_1(\Cp2\sminus B')=\RBG3$.
\qed
\endcorollary

\lemma\label{torus.epi}
Let~$B$ be a sextic of torus type with
$\pi_1(\Cp2\sminus B)=\RBG3$,
and let~$\MB$ be a Milnor ball about an inner singular point~$P$
of~$B$. Then the inclusion homomorphism
$\pi_1(\MB\sminus B)\to\pi_1(\Cp2\sminus B)$ is onto.
\endlemma

\proof
Consider a perturbation~$B_t$, $t\in[0,1]$, of $B=B_0$ to a six
cuspidal sextic $B'=B_1$, transversal to $\partial\MB$.
%Let $P'\in\MB$ be a cusp of~$B'$, let~$\MB'$ be a Milnor ball
%about~$P'$,
Pick a cusp $P'\in\MB$ of~$B'$, consider a Milnor ball~$\MB'$
about~$P'$,
and let $i\:\MB'\sminus B'\into\MB\sminus B'$ and
$j\:\MB\sminus B'\into\Cp2\sminus B'$ be the inclusions.
The homomorphism
$(j\circ i)_*\:\pi_1(\MB'\sminus B')\to\pi_1(\Cp2\sminus B')$ is
onto, see~\cite{Zariski.group};
hence, so is $j_*$. On the other
hand, the perturbation $B\to B'$ induces an epimorphism
$\pi_1(\MB\sminus B)\onto\pi_1(\MB\sminus B')$ and an isomorphism
$\pi_1(\Cp2\sminus B)=\pi_1(\Cp2\sminus B')$, see
Corollary~\ref{torus->torus}; hence,
the inclusion homomorphism
$\pi_1(\MB\sminus B)\to\pi_1(\Cp2\sminus B)$ is onto.
\endproof

\lemma\label{torus->nontorus}
Let~$B$ be a sextic of torus type with simple singularities and
with $\pi_1(\Cp2\sminus B)=\RBG3$,
and let~$B'$ be a perturbation of~$B$ which is not of torus type.
Then $\pi_1(\Cp2\sminus B')=\CG6$.
\endlemma

\proof
Since $B'$ is not of torus type, there is an inner singular
point~$P$ of~$B$ that undergoes a perturbation not of torus type.
Let~$\MB$ be a Milnor ball about~$P$. If $P$ is of type~$\bE_6$,
the group $\pi_1(\MB\sminus B')$ is abelian, see~\cite{degt.e6},
and the statement follows from Lemma~\ref{torus.epi}.
Assume that $P$ is of type $\bA_{3k-1}$. Then
$$
\pi_1(\MB\sminus B)=\<\Ga,\Gb\,|\,\Gs^{3k}=\id>
\quad\text{and}\quad
\pi_1(\MB\sminus B')=\<\Ga,\Gb\,|\,\Gs^s=\id>
$$
for some $s\ne0\bmod3$, see the proof of
Lemmas~\ref{pert.a11}--\ref{pert.a7}. On the other hand, in the
group $\pi_1(\Cp2\sminus B)=\RBG3$, the generators~$\Ga$, $\Gb$
are subject to the braid relation $\Gs^3=\id$. (From the proof of
Lemma~\ref{torus.epi}, it follows that $\Ga$, $\Gb$ are taken to
the standard generators $\Gs_1,\Gs_2\in\RBG3$.) Hence, after the
perturbation, one has a relation $\Gs=\id$,
\ie, $\Ga=\Gb$ in $\pi_1(\Cp2\sminus B')$, and
Lemma~\ref{torus.epi} applies.
\endproof

\proof[Proof of Theorem~\ref{th.group}]
Due to the results of Sections~\ref{s.a17+a2}
and~\ref{s.2a8+a3}, Corollary~\ref{torus->torus} implies that, for any
sextic~$B$ of torus type obtained by a perturbation from an
irreducible sextic with the set of singularities
$(\bA_{17})\splus\bA_2$ or $(2\bA_8)\splus\bA_3$, one has
$\pi_1(\Cp2\sminus B)=\RBG3$. In particular, this statement covers
all curves listed in Theorem~\ref{th.group}, as their
moduli spaces are connected, see Theorem~\ref{th.moduli}.
\endproof

Now, consider a perturbation~$B'$ of~$B$ that is \emph{not} of
torus type. The extremal sets of singularities obtained in this way
are
$$
\gather
\bA_{16}\splus\bA_2,\quad
 \bA_{15}\splus\bA_2\splus\bA_1,\quad
 \bA_{13}\splus\bA_3\splus\bA_2,\\\allowdisplaybreak
 \bA_{12}\splus\bA_4\splus\bA_2,\quad
 \bA_{10}\splus\bA_6\splus\bA_2,\quad
 \bA_9\splus\bA_7\splus\bA_2,\\\allowdisplaybreak
\bA_8\splus\bA_7\splus\bA_3,\quad
 \bA_8\splus\bA_6\splus\bA_3\splus\bA_1,\quad
 \bA_8\splus\bA_4\splus2\bA_3.
\endgather
$$
Due to Lemma~\ref{torus->nontorus}, each of these sets of
singularities is realized by a sextic whose fundamental
group is cyclic.

\section{Reducible curves of torus type}\label{S.groups}

Now, we consider the maximal reducible sextics of torus type of
the form $B=\DC(\B_2,\L)$, where $\B_2$ is the trigonal curve as
in~\ref{s.a5+a2}.
The computation of the group $\pi_1(\Cp2\sminus B)$
follows the outline in Sections~\ref{s.group}
and~\ref{s.monodromy}, with an additional
simplification due to the fact that the singular
fibers of $\B_2+\L$ are all real;
we systematically ignore the relation from the singular fiber at
infinity.

Instead of simplifying the obtained presentations for
$\pi_1(\Cp2\sminus B)$, we perturb~$B$ to an irreducible
sextic~$B'$ and use Lemmas~\ref{pert.a11}--\ref{pert.e7} to
compute $\pi_1(\Cp2\sminus B')$.
If $B'$ is of torus type, we only prove that there is an
epimorphism $\BG3\onto\pi_1(\Cp2\sminus B)$;
Lemma~\ref{RBG->torus} implies that it is an isomorphism.
(Similarly, if $B'$ is not of torus type, we only prove that the
group is abelian.)
Furthermore, we only consider maximal irreducible perturbations;
as the groups obtained are $\RBG3$ or $\CG6$, the results extend
to other perturbations using Corollary~\ref{torus->torus} and
Lemma~\ref{torus->nontorus}.

Here and in~\S\ref{S.others}, without further references, the
perturbations are constructed using Proposition~5.1.1
in~\cite{degt.8a2}, by perturbing the singular points
independently.

\midinsert
\plot{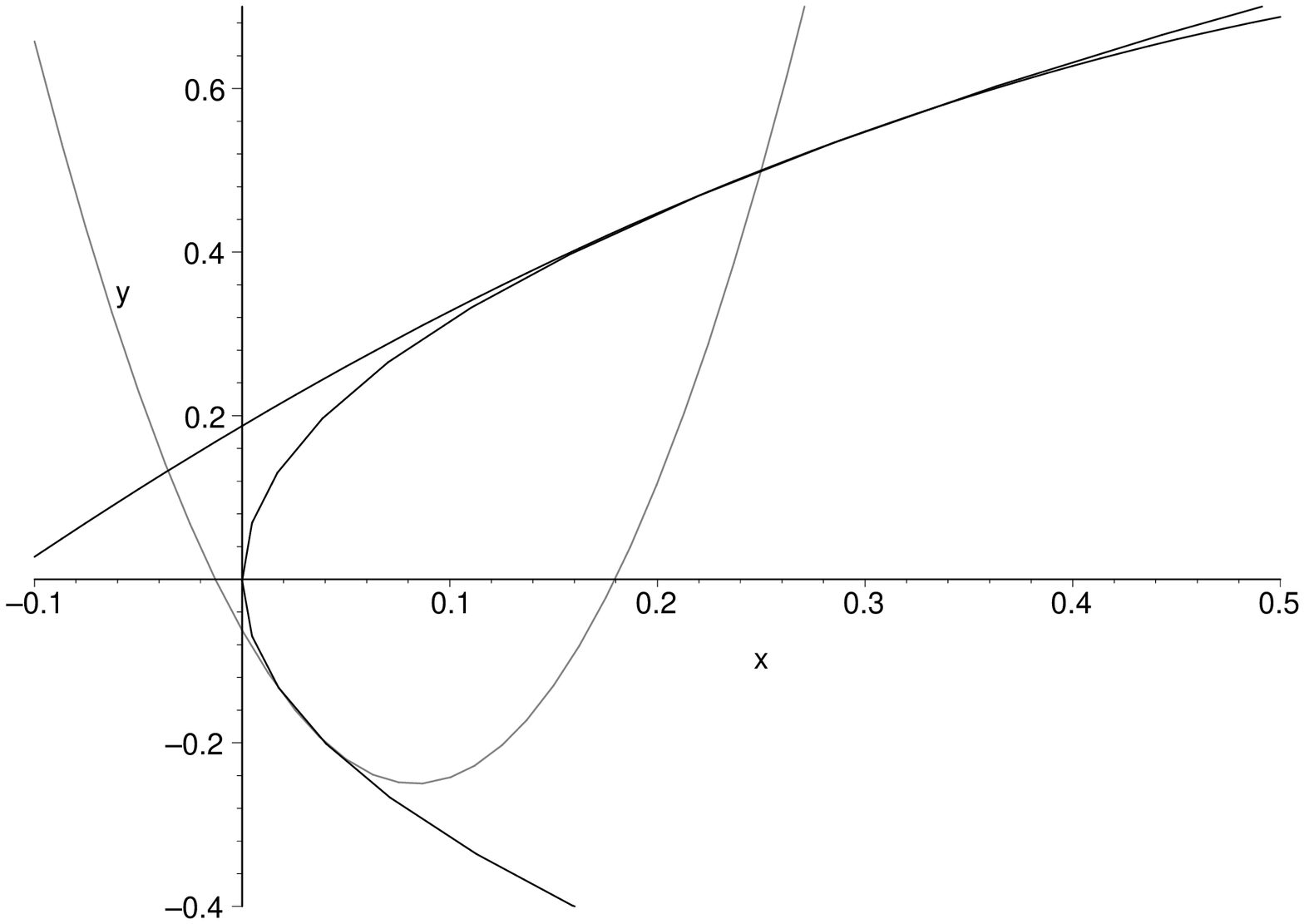}
\figure
The set of singularities
$(\bA_{11}\splus2\bA_2)\splus\bA_2\splus2\bA_1$
\endfigure\label{fig.a11+3a2}
\endinsert

\subsection{The set of singularities
$(\bA_{11}\splus2\bA_2)\splus\bA_2\splus2\bA_1$}\label{s.a11+3a2}
Take for~$\L$ the section given by~\eqref{eq.a11+3a2}, see
Figure~\ref{fig.a11+3a2} (where the point~$R_1$ of transversal
intersection of~$\L'$ and the lower branch of~$\B_2'$
is missing), and
choose the generators
%$\Ga$, $\Gb$, $\Gd$, $\Gg$
$(\eta_1,\eta_2,\eta_3,\eta_4)=(\Ga,\Gb,\Gd,\Gg)$
in a real fiber~$F$ just to the left from~$R_5$ (\eg, over $x=0.2$).
The relations for $\bpi$ are:
$$
\alignat2
&\Gd(\Ga\Gb)^3=(\Gb\Ga\Gb)\Gd(\Ga\Gb\Ga)&\qquad&
 \text{(the fiber through~$R_5$)},\\\allowdisplaybreak
&[\Gd,\Ga\Gb]=1&&
 \text{(the fiber through~$R_5$)},\\\allowdisplaybreak
&[(\Gb\Ga\Gb\Gd)\1\Ga(\Gb\Ga\Gb\Gd),\Gg]=1&&
 \text{(the fiber through~$R_1$)},\\\allowdisplaybreak
&(\Gg\Gd)^3=(\Gd\Gg)^3&&
 \text{(the fiber through~$Q_5$)},\\\allowdisplaybreak
&\Gb=(\Gd\Gg\Gd)\Gg(\Gd\Gg\Gd)\1&&
 \text{(the vertical tangent)},\\\allowdisplaybreak
&[\Gb\1\Ga\Gb,(\Gd\Gg)\Gd(\Gd\Gg)\1]=1&&
 \text{(the fiber through~$Q_1$)},\\\allowdisplaybreak
&(\Ga\Gb\Gd\Gg)^2=1&&
 \text{(the relation at infinity)}.
\endalignat
$$
Passing to the generators
$\Ga$, $\bGa$, $\Gb$, $\bGb$, $\Gg$, $\bGg$, see
Lemma~\ref{bpi->pi}, we obtain the following set of relations for
$\pi_1(\Cp2\sminus B)$, where $B$ is a reducible sextic of torus
type with the set of singularities
$(\bA_{11}\splus2\bA_2)\splus\bA_2\splus2\bA_1$:
$$
\gather
(\Ga\Gb)^3=\bGb\bGa\bGb\Ga\Gb\Ga,\quad
 (\bGa\bGb)^3=\Gb\Ga\Gb\bGa\bGb\bGa,
 \eqtag\label{eq.4.1}\\\allowdisplaybreak
\Ga\Gb=\bGa\bGb,
 \eqtag\label{eq.4.2}\\\allowdisplaybreak
[(\Gb\Ga\Gb)\1\Ga(\Gb\Ga\Gb),\bGg]=1,\quad
 [(\bGb\bGa\bGb)\1\bGa(\bGb\bGa\bGb),\Gg]=1,
 \eqtag\label{eq.4.3}\\\allowdisplaybreak
\Gg\bGg\Gg=\bGg\Gg\bGg,
 \eqtag\label{eq.4.4}\\\allowdisplaybreak
\Gb=\bGg\Gg\bGg\1,\quad
 \bGb=\Gg\bGg\Gg\1
 \eqtag\label{eq.4.5}\\\allowdisplaybreak
(\bGb\Gg)\1\bGa\bGb\Gg=(\Gb\bGg)\1\Ga\Gb\bGg,
 \eqtag\label{eq.4.6}\\\allowdisplaybreak
\Ga\Gb\bGg\bGa\bGb\Gg=1.
 \eqtag\label{eq.4.7}
\endgather
$$
In the presence of~\eqref{eq.4.1} and~\eqref{eq.4.2},
relations~\eqref{eq.4.3} simplify to
$$
[(\bGa\bGb)\bGa(\bGa\bGb)\1,\bGg]=1,\quad
 [(\Ga\Gb)\Ga(\Ga\Gb)\1,\Gg]=1.
 \eqtag\label{eq.4.8}
$$

Consider the perturbation $\bA_{11}\mapsto\bA_8\splus\bA_2$,
producing an irreducible sextic~$B'$ of torus type with the set of
singularities
$$
(\bA_8\splus3\bA_2)\splus\bA_2\splus2\bA_1.
$$
%and, after a further perturbation of the outer singular points,
%also
%$$
%(\bA_8\splus3\bA_2)\splus\bA_2\splus\bA_1
%\quad\text{and}\quad
%(\bA_8\splus3\bA_2)\splus3\bA_1.
%$$
This perturbation adds the braid relation
$\Ga\Gb\Ga=\Gb\Ga\Gb$. Then, the first relation in~\eqref{eq.4.1}
simplifies to $\Gb\Ga\Gb=\bGb\bGa\bGb$; in view of~\eqref{eq.4.2},
this implies that $\Gb=\bGb$. Similarly, using the second
relation in~\eqref{eq.4.1}, one obtains $\Ga=\bGa$. Furthermore,
in the presence of the braid relation,
one has $(\Ga\Gb)\Ga(\Ga\Gb)\1=\Gb$;
hence, \eqref{eq.4.8}
implies $[\Gb,\bGg]=[\Gb,\Gg]=1$ and \eqref{eq.4.5} yields
$\Gg=\bGg=\Gb$. Thus, the map $\Ga,\bGa\mapsto\Gs_1$,
$\Gb,\bGb,\Gg,\bGg\mapsto\Gs_2$ establishes an isomorphism
$\pi_1(\Cp2\sminus B')=\RBG3$.

Now, perturb one of the nodes over~$R_1$,
producing an irreducible sextic~$B'$ of
torus type with the set of singularities
$$
(\bA_{11}\splus2\bA_2)\splus\bA_2\splus\bA_1.
%\quad\text{and}\quad
%(\bA_{11}\splus2\bA_2)\splus2\bA_1.
$$
%(The latter set of singularities is obtained by a further
%perturbation of the outer cusp to a node.)
This perturbation
simplifies one of the two relations~\eqref{eq.4.8}: for example,
we can replace the first one with
$$
\bGg=(\bGa\bGb)\bGa(\bGa\bGb)\1.
\eqtag\label{eq.4.9}
$$
To simplify the group, introduce the generators $u=\Ga\Gb$ and
$v=\Ga\Gb\Ga$, so that $\Ga=u\1v$ and $\Gb=v\1u^2$. Then $\bu=u$
(from~\eqref{eq.4.2}) and $\bv=u^{-3}vu^3$ (from the second
relation in~\eqref{eq.4.1}). Hence, $\bGa=u^{-4}vu^3$,
$\bGb=u^{-3}v\1u^5$, and the new relation~\eqref{eq.4.9} turns
into $\bGg=u^{-3}vu^2$; in view of~\eqref{eq.4.7}, this implies
$\Gg=u^{-3}v\1u^2$. Substituting the expressions obtained to the
second relation in~\eqref{eq.4.5}, we arrive at $u^3=v^2$; hence,
$\bGa=\Ga=u\1v$, $\bGb=\Gb=\bGg=v\1u^2$, and $\Gg=v^{-3}u^2$.
Substituting to the first relation in~\eqref{eq.4.5}, we get
$v^2=1$. Thus, also $u^3=1$, and we obtain an isomorphism
$\pi_1(\Cp2\sminus B')=\RBG3$, given by
$\Ga,\bGa\mapsto\Gs_1$ and $\Gb,\bGb,\Gg,\bGg\mapsto\Gs_2$.

Finally, consider a maximal irreducible perturbation of the
type~$\bA_{11}$ singular point that is not of torus type,
producing irreducible plane sextics with the sets of singularities
$$
\bA_{10}\splus3\bA_2\splus2\bA_1,\quad
 \bA_6\splus\bA_4\splus3\bA_2\splus2\bA_1
$$
(see Lemma~\ref{pert.a11}). According to Lemma~\ref{pert.a11},
this perturbation introduces the additional relations
$\Ga=\bGa=\Gb=\bGb$. Then, from~\eqref{eq.4.3} one has
$[\Ga,\Gg]=[\Ga,\bGg]=1$ and, due to~\eqref{eq.4.7}, also
$[\Gg,\bGg]=1$. Hence, the group is abelian.

\subsection{The set of singularities
$(2\bA_5\splus2\bA_2)\splus\bD_5$}\label{s.2a5+2a2+d5}
Consider the triple $\B_2'$, $\L'$, $\L$ as in
Section~\ref{s.a11+3a2}, see Figure~\ref{fig.a11+3a2}, but now let
$\B=\B_2'+\L$ and
%take~$\L'$ for the section defining the double
%covering $\Cp2\to\Sigma_2/E$.
$B=\DC((\B_2'+\L),\L')$: it is
a reducible sextic
of torus type with the set of singularities
$(2\bA_5\splus2\bA_2)\splus\bD_5$.
The group
$\bpi=\pi_1(\Sigma_2\sminus(\B\cup\L'\cup E))$ is found in
Section~\ref{s.a11+3a2}. Modifying Lemma~\ref{bpi->pi},
let $\Ga^2=1$ and pass to the
generators $\Gb$, $\bGb=\Ga\Gb\Ga$, $\Gg$, $\bGg=\Ga\Gg\Ga$,
$\Gd$, $\bGd=\Ga\Gd\Ga$. We obtain the following set of relations for
the group $\pi_1(\Cp2\sminus B)$:
$$
\gather
\bGd\Gb\bGb\Gb=\Gd\bGb\Gb\bGb,
 \eqtag\label{eq.3.1}\\\allowdisplaybreak
\bGd\Gb=\Gb\Gd,\quad
 \Gd\bGb=\bGb\bGd,
 \eqtag\label{eq.3.2}\\\allowdisplaybreak
(\Gb\bGb\bGd)\bGg(\Gb\bGb\bGd)\1=(\bGb\Gb\Gd)\Gg(\bGb\Gb\Gd)\1,
 \eqtag\label{eq.3.3}\\\allowdisplaybreak
(\Gg\Gd)^3=(\Gd\Gg)^3,\quad
 (\bGg\bGd)^3=(\bGd\bGg)^3,
 \eqtag\label{eq.3.4}\\\allowdisplaybreak
\Gb=(\Gd\Gg\Gd)\Gg(\Gd\Gg\Gd)\1,\quad
 \bGb=(\bGd\bGg\bGd)\bGg(\bGd\bGg\bGd)\1,
 \eqtag\label{eq.3.5}\\\allowdisplaybreak
(\Gb\Gd\Gg)\Gd(\Gb\Gd\Gg)\1=(\bGb\bGd\bGg)\bGd(\bGb\bGd\bGg)\1,
 \eqtag\label{eq.3.6}\\\allowdisplaybreak
\bGb\bGd\bGg\Gb\Gd\Gg=1.
 \eqtag\label{eq.3.7}
\endgather
$$
(Here, \eqref{eq.3.1} is simplified using~\eqref{eq.3.2}.)

The perturbation $\bD_5\mapsto\bA_4$ produces an
irreducible sextic~$B'$ of torus type with the set of singularities
$$
(2\bA_5\splus2\bA_2)\splus\bA_4
$$
and introduces the relation $\Gb=\bGb=\Gd=\bGd$, see
Lemma~\ref{pert.d5}. Then, due to~\eqref{eq.3.3}, $\Gg=\bGg$ and
\eqref{eq.3.5} simplifies to $\Gb\Gg\Gb=\Gg\Gb\Gg$.
Hence,
%$\pi_1(\Cp2\sminus B')=\RBG3$, the
%isomorphism being $\Gb,\bGb,\Gd,\bGd\mapsto\Gs_1$,
%$\Gg,\bGg\mapsto\Gs_2$.
the map $\Gb,\bGb,\Gd,\bGd\mapsto\Gs_1$,
$\Gg,\bGg\mapsto\Gs_2$
establishes an isomorphism $\pi_1(\Cp2\sminus B')=\RBG3$.

The perturbation $\bA_5\mapsto\bA_4$ produces an
irreducible sextic~$B'$
%, not of torus type,
with the set of singularities
$$
\bD_5\splus\bA_5\splus\bA_4\splus2\bA_2
$$
and introduces the relation $\Gd=\Gg$, see Lemma~\ref{pert.a5}.
Then, due to the first relation in~\eqref{eq.3.5}, $\Gb=\Gg$,
relation~\eqref{eq.3.2} implies $\bGd=\Gg$ and $[\bGb,\Gg]=1$, and
one has $\bGb=\Gg$ (from~\eqref{eq.3.1}) and $\bGg=\Gg^{-5}$
(from~\eqref{eq.3.7}). Thus, the group is abelian.

\midinsert
\plot{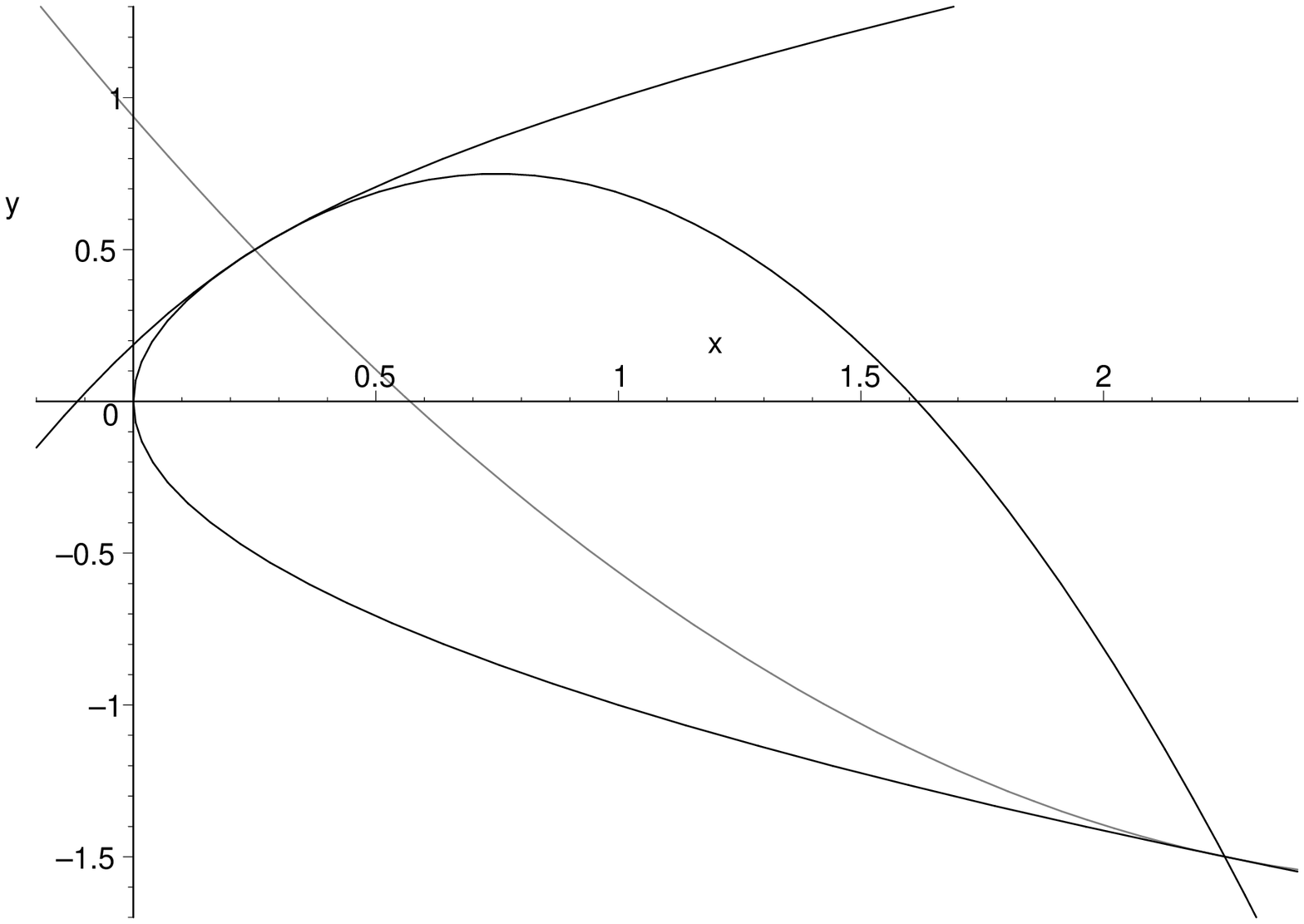}
\figure
The set of singularities
$(\bA_{11}\splus2\bA_2)\splus\bD_4$
\endfigure\label{fig.a11+2a2+d4}
\endinsert

\subsection{The set of singularities
$(\bA_{11}\splus2\bA_2)\splus\bD_4$}\label{s.a11+2a2+d4}
Take for~$\L$ the section given by~\eqref{eq.a11+2a2+d4}, see
Figure~\ref{fig.a11+2a2+d4} (where the point~$Q_1$ of transversal
intersection of~$\L$ and the upper branch of~$\B_2'$ is missing).
Choose the generators
$(\eta_1,\eta_2,\eta_3,\eta_4)=(\Ga,\Gb,\Gd,\Gg)$
%$\Ga$, $\Gb$, $\Gd$, $\Gg$
in a real fiber~$F$ between~$R_5$ and~$R_1$ (\eg, over $x=1$).
The relations are:
$$
\alignat2
&\Gd(\Ga\Gb)^3=(\Gb\Ga\Gb)\Gd(\Ga\Gb\Ga)&\qquad&
 \text{(the fiber through~$R_5$)},\\\allowdisplaybreak
&[\Gd,\Ga\Gb]=1&&
 \text{(the fiber through~$R_5$)},\\\allowdisplaybreak
&\Gb(\Gd\Gg)^2=\Gg\Gb\Gd\Gg\Gd&&
 \text{(the fiber through~$R_1$)},\\\allowdisplaybreak
&[\Gb,\Gd\Gg]=1&&
 \text{(the fiber through~$R_1$)},\\\allowdisplaybreak
&(\Gb\Ga\Gb\Gd)\1\Ga(\Gb\Ga\Gb\Gd)=\Gg&&
 \text{(the vertical tangent)},\\\allowdisplaybreak
&[\Ga,\Gg\1\Gd\Gg]=1&&
 \text{(the fiber through~$Q_1$)},\\\allowdisplaybreak
&(\Ga\Gb\Gd\Gg)^2=1&&
 \text{(the relation at infinity)}.
\endalignat
$$
Passing to the generators
$\Ga$, $\bGa$, $\Gb$, $\bGb$, $\Gg$, $\bGg$, see
Lemma~\ref{bpi->pi}, we obtain the
following relations for the group $\pi_1(\Cp2\sminus B)$ of a
reducible sextic~$B$ of torus type with the set of singularities
$(\bA_{11}\splus2\bA_2)\splus\bD_4$:
$$
\gather
(\Ga\Gb)^3=\bGb\bGa\bGb\Ga\Gb\Ga,\quad
 (\bGa\bGb)^3=\Gb\Ga\Gb\bGa\bGb\bGa,
 \eqtag\label{eq.5.1}\\\allowdisplaybreak
\Ga\Gb=\bGa\bGb,
 \eqtag\label{eq.5.2}\\\allowdisplaybreak
\Gb\bGg\Gg=\Gg\Gb\bGg,\quad
 \bGb\Gg\bGg=\bGg\bGb\Gg,
 \eqtag\label{eq.5.3}\\\allowdisplaybreak
\Gb\bGg=\bGg\bGb,\quad
 \bGb\Gg=\Gg\Gb,
 \eqtag\label{eq.5.4}\\\allowdisplaybreak
(\Gb\Ga\Gb)\1\Ga(\Gb\Ga\Gb)=\bGg,\quad
 (\bGb\bGa\bGb)\1\bGa(\bGb\bGa\bGb)=\Gg,
 \eqtag\label{eq.5.5}\\\allowdisplaybreak
\Ga\Gg\1\bGg=\Gg\1\bGg\bGa,\quad
 \bGa\bGg\1\Gg=\bGg\1\Gg\Ga,
 \eqtag\label{eq.5.6}\\\allowdisplaybreak
\Ga\Gb\bGg\bGa\bGb\Gg=1.
 \eqtag\label{eq.5.7}
\endgather
$$
As above, in view of~\eqref{eq.5.1} and~\eqref{eq.5.2},
relations~\eqref{eq.5.5} simplify to
$$
(\bGa\bGb)\bGa(\bGa\bGb)\1=\bGg,\quad
 (\Ga\Gb)\Ga(\Ga\Gb)\1=\Gg.
 \eqtag\label{eq.5.8}
$$

The perturbation $\bA_{11}\mapsto\bA_8\splus\bA_2$
produces an irreducible sextic~$B'$ of torus type with the set of
singularities
$$
(\bA_8\splus3\bA_2)\splus\bD_4.
$$
The perturbation adds to the presentation a braid relation
$\Ga\Gb\Ga=\Gb\Ga\Gb$. Similar to Section~\ref{s.a11+3a2} (the
perturbation $\bA_{11}\to\bA_8\splus\bA_2$),
relations~\eqref{eq.5.1} and~\eqref{eq.5.2} imply that $\Ga=\bGa$ and
$\Gb=\bGb$, and~\eqref{eq.5.8} turns into $\Gg=\bGg=\Gb$. Hence,
there is an isomorphism $\pi_1(\Cp2\sminus B')=\RBG3$
given by the map $\Ga,\bGa\mapsto\Gs_1$,
$\Gb,\bGb,\Gg,\bGg\mapsto\Gs_2$.

Now, perturb the type~$\bD_4$ point of~$B$.
After the perturbation,
the generators $\Gb$, $\bGb$, $\Gg$, and~$\bGg$ pairwise commute,
see Lemma~\ref{pert.d4}.
It follows that $\Gb=\bGb$ (from~\eqref{eq.5.4}), $\Ga=\bGa$
(from~\eqref{eq.5.2}), and $\Gg=\bGg$ (from~\eqref{eq.5.8}), and
the presentation simplifies to
$$
\bigl<\Ga,\Gb\bigm|
 (\Ga\Gb)^3=(\Gb\Ga)^3,\ (\Ga\Gb\Ga)^2=1\bigr>.
\eqtag\label{eq.D4.perturbed}
$$
This is, indeed, the fundamental group of a reducible curve of
torus type whose set of singularities is
%$(\bA_{11}\splus2\bA_2)\splus3\bA_1$,
%$(\bA_{11}\splus2\bA_2)\splus\bA_3$, or
%$(\bA_{11}\splus2\bA_2)\splus2\bA_1$.
$$
(\bA_{11}\splus2\bA_2)\splus3\bA_1,\quad
(\bA_{11}\splus2\bA_2)\splus\bA_3,\quad\text{or}\quad
(\bA_{11}\splus2\bA_2)\splus2\bA_1.
$$
To obtain an irreducible curve with the set of singularities
$$
(\bA_{11}\splus2\bA_2)\splus\bA_3,
$$
we choose the perturbation
$\bD_4\mapsto\bA_3$ so that the generators $\Gb$ and~$\Gg$ become
conjugate (hence equal) in $\pi_1(\MB_P\sminus B')$. (Locally, we
perturb a triple of lines to a conic and a line tangent to it, and
the choice of a line to be kept as a separate component is
governed by the choice of a subdiagram $\bA_3\subset\bD_4$; any
such choice can be realized by a perturbation of~$B$.)
Then,
\eqref{eq.5.8} implies the braid relation $\Ga\Gb\Ga=\Gb\Ga\Gb$,
and the map $\Ga,\bGa\mapsto\Gs_1$,
$\Gb,\bGb,\Gg,\bGg\mapsto\Gs_2$ establishes an isomorphism
$\pi_1(\Cp2\sminus B')=\RBG3$.

The perturbations $\bA_{11}\mapsto\bA_{10}$ or $\bA_6\splus\bA_4$,
see Lemma~\ref{pert.a11}, produce irreducible sextics with the
sets of singularities
$$
\bD_4\splus\bA_{10}\splus2\bA_2,\quad
 \bD_4\splus\bA_6\splus\bA_4\splus2\bA_2
$$
while adding the relation $\Ga=\bGa=\Gb=\bGb$. Then \eqref{eq.5.5}
implies $\Gg=\bGg=\Ga$, and the group is abelian.

%\midinsert
%\plot{a11+e6}
%\figure
%The set of singularities
%$(\bE_6\splus\bA_{11})\splus2\bA_1$
%\endfigure\label{fig.a11+e6}
%\endinsert

\midinsert
\plot{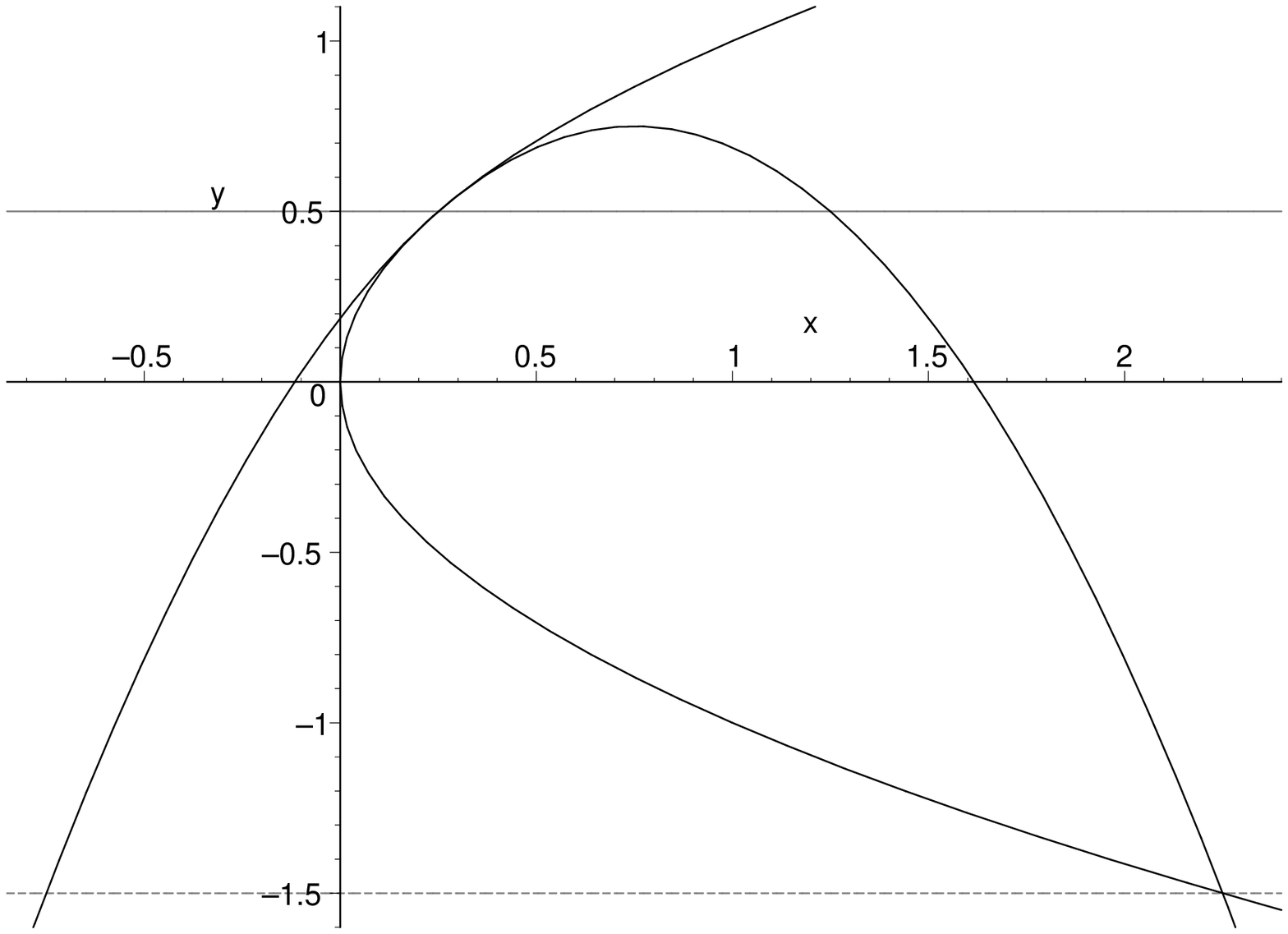}
\figure
The sets
%of singularities
$(\bE_6\splus\bA_{11})\splus2\bA_1$ and
$(\bE_6\splus2\bA_5)\splus\bA_3$
\endfigure\label{fig.horizontal}
\endinsert

\subsection{The set of singularities
$(\bE_6\splus\bA_{11})\splus2\bA_1$}\label{s.a11+e6}
Let~$\L$ be the section given by~\eqref{eq.a11+e6}, see the solid
horizontal grey line in Figure~\ref{fig.horizontal}. (The resulting
set of singularities $(\bE_6\splus\bA_{11})\splus2\bA_1$
is erroneously missing in Oka~\cite{Oka.reducible}.)
Choosing the generators
%$\Ga$, $\Gb$, $\Gd$, $\Gg$
$(\eta_1,\eta_2,\eta_3,\eta_4)=(\Ga,\Gb,\Gd,\Gg)$
in a real fiber~$F$ between~$R_5$ and~$Q_1$ (\eg, over $x=1$),
we obtain the relations
$$
\alignat2
&[\Gd,\Gb]=1&\qquad&
 \text{(the fiber through~$Q_1$)},\\\allowdisplaybreak
&[\Gd,\Ga\Gb]=1&&
 \text{(the fiber through~$R_5$)}.
\endalignat
$$
(We only list the few relations needed in the sequel.) Hence, also
$[\Gd,\Ga]=1$. The relation at the vertical
tangent of~$\B$ has the form $\Gg=\text{(a word in $\Ga$, $\Gb$,
$\Gd$)}$; hence, we also have $[\Gd,\Gg]=1$. Thus, $\Gd$ is a
central element, and a presentation for
the group
$\pi_1(\Sigma_2\sminus(\B\cup\L\cup E))$ can be obtained from
%that
the presentation
in Section~\ref{s.a11+2a2+d4} by adding the relations
$[\Ga,\Gd]=[\Gb,\Gd]=[\Gg,\Gd]=1$. After eliminating
$\Gg=(\Gb\Ga)\1\Ga(\Gb\Ga)$, we get
$$
\bigl<\Ga,\Gb,\Gd\bigm|
 (\Ga\Gb)^3=(\Gb\Ga)^3,\
 [\Ga,\Gd]=[\Gb,\Gd]=1,\
 (\Ga\Gb\Ga)^2\Gd^2=1\bigr>.
\eqtag\label{eq.2E7.1}
$$
The fundamental group $\pi_1(\Cp2\sminus B)$ of a reducible
sextic~$B$ of torus type with the set of singularities
$(\bE_6\splus\bA_{11})\splus2\bA_1$
is obtained from~\eqref{eq.2E7.1}
by letting $\Gd=1$. The result
is~\eqref{eq.D4.perturbed}.

A perturbation of a node of~$B$ produces an irreducible sextic
of torus type with the set of singularities
$$
(\bE_6\splus\bA_{11})\splus\bA_1.
$$
The additional relation introduced by this operation is $\Gb=\Gg$.
Hence, the resulting fundamental group is $\RBG3$,
\cf. the perturbation $\bD_4\mapsto\bA_3$ in
Section~\ref{s.a11+2a2+d4}.

The perturbation $\bA_{11}\mapsto\bA_8\splus\bA_2$ produces
the set of singularities
$$
(\bE_6\splus\bA_8\splus\bA_2)\splus2\bA_1.
$$
The additional relation is $\Ga\Gb\Ga=\Gb\Ga\Gb$, and the
resulting group is $\RBG3$.

The perturbations $\bA_{11}\mapsto\bA_{10}$ or $\bA_6\splus\bA_4$
produce the sets of singularities
$$
\bE_6\splus\bA_{10}\splus2\bA_1,\quad
 \bE_6\splus\bA_6\splus\bA_4\splus2\bA_1,
$$
while adding to~\eqref{eq.D4.perturbed} the relation $\Ga=\Gb$, see
Lemma~\ref{pert.a11}. Hence, the resulting group is abelian.

%\midinsert
%\plot{2a5+e6+a3}
%\figure
%The set of singularities
%$(\bE_6\splus2\bA_5)\splus\bA_3$
%\endfigure\label{fig.2a5+e6+a3}
%\endinsert

\subsection{The set of singularities
$(\bE_6\splus2\bA_5)\splus\bA_3$}\label{s.2a5+e6+a3}
Let~$\L$ be the section given by~\eqref{eq.2a5+e6+a3}, see the
dotted horizontal grey line in Figure~\ref{fig.horizontal}.
Choosing the generators
$(\eta_1,\eta_2,\eta_3,\eta_4)=(\Ga,\Gb,\Gg,\Gd)$
%$\Ga$, $\Gb$, $\Gg$, $\Gd$
in a real fiber~$F$ between~$R_5$ and~$R_1$ (\eg, over $x=1$), we
obtain the relations
$$
\alignat2
&(\Ga\Gb)^3=(\Gb\Ga)^3&\qquad&
 \text{(the fiber through~$R_5$)},\\\allowdisplaybreak
&(\Ga\Gb)\Ga(\Ga\Gb)\1=\Gg&&
 \text{(the vertical tangent)},\\\allowdisplaybreak
&[\Gd,\Ga\Gb\Ga\1]=1&&
 \text{(the fiber through~$Q_1$)},\\\allowdisplaybreak
&[\Gg\Gd,\Gb]=[\Gb\Gg\Gd,\Gg]=[\Gb\Gg,\Gd]=1&&
 \text{(the fiber through~$R_1$)},\\\allowdisplaybreak
&(\Ga\Gb\Gg\Gd)^2=1&&
 \text{(the relation at infinity)}.
\endalignat
$$
(The second and the third relations were simplified using the
first one.)
Passing to $\Ga$, $\bGa$, $\Gb$, $\bGb$, $\Gg$, $\bGg$, see
Lemma~\ref{bpi->pi}, gives the following relations for the group
of a reducible sextic~$B$ of torus type with the set of
singularities $(\bE_6\splus2\bA_5)\splus\bA_3$:
$$
\gather
(\Ga\Gb)^3=(\Gb\Ga)^3,\quad
 (\bGa\bGb)^3=(\bGb\bGa)^3,
 \eqtag\label{eq.6.1}\\\allowdisplaybreak
(\Ga\Gb)\Ga(\Ga\Gb)\1=\Gg,\quad
 (\bGa\bGb)\bGa(\bGa\bGb)\1=\bGg,
 \eqtag\label{eq.6.2}\\\allowdisplaybreak
\Ga\Gb\Ga\1=\bGa\bGb\bGa\1,
 \eqtag\label{eq.6.3}\\\allowdisplaybreak
\Gg\bGb=\Gb\Gg=\bGb\bGg=\bGg\Gb,
 \eqtag\label{eq.6.4}\\\allowdisplaybreak
\Ga\Gb\Gg\bGa\bGb\bGg=1.
 \eqtag\label{eq.6.5}
\endgather
$$
The perturbation $\bA_5\mapsto\bA_4$ of one of the type~$\bA_5$
singular points of~$B$ produces the set of singularities
$$
\bE_6\splus\bA_5\splus\bA_4\splus\bA_3
$$
and adds the relation $\Ga=\Gb$, see Lemma~\ref{pert.a5}. Then
$\Gg=\Ga$ (from~\eqref{eq.6.2}), $\bGb=\bGg=\Ga$
(from~\eqref{eq.6.4}), and $\bGa=\Ga^{-5}$ (from~\eqref{eq.6.5}).
Hence, the group is abelian.

%\midinsert
%\plot{3a5+d4}
%\figure
%The set of singularities
%$(3\bA_5)\splus\bD_4$
%\endfigure\label{fig.3a5+d4}
%\endinsert

\midinsert
\plot{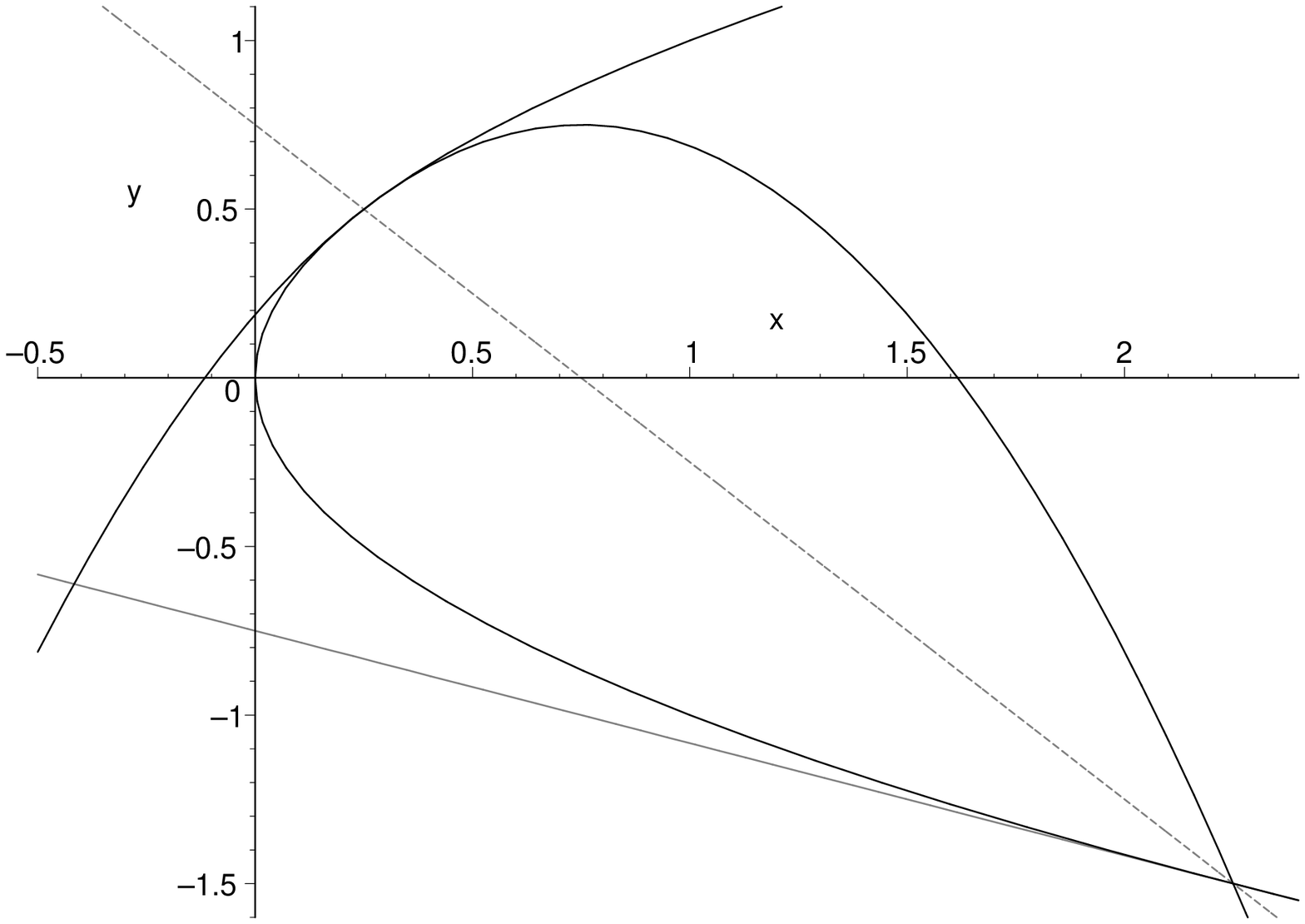}
\figure
The sets
%of singularities
$(3\bA_5)\splus\bD_4$ and
$(\bA_{11}\splus\bA_5)\splus\bA_3$
\endfigure\label{fig.slanted}
\endinsert

\subsection{The set of singularities
$(3\bA_5)\splus\bD_4$}\label{s.3a5+d4}
Let~$\L$ be the section given by~\eqref{eq.3a5+d4}, see the
solid grey line in Figure~\ref{fig.slanted}.
Choosing the generators $\Ga$, $\Gb$, $\Gg$, $\Gd$ as in
Section~\ref{s.2a5+e6+a3}, we obtain the same set of relations,
except that the relation at the fiber through~$R_1$ should be
replaced with
$$
[\Gb,\Gg\Gd]=1,\quad
 \Gd\Gb\Gg\Gd\Gg=\Gb\Gg\Gd\Gg\Gd\quad
 \text{(the fiber through~$R_1$)}.
$$
Hence, the relations for $\pi_1(\Cp2\sminus B)$ are
\eqref{eq.6.1}--\eqref{eq.6.3}, \eqref{eq.6.5}, and the relations
$$
\Gb\Gg=\Gg\bGb,\quad
 \bGb\bGg=\bGg\Gb,\quad
 \bGb\bGg\Gg=\Gb\Gg\bGg
 \eqtag\label{eq.7.1}
$$
replacing~\eqref{eq.6.4}.

We are only interested in the perturbation $\bA_5\mapsto\bA_4$ of
two of the three type~$\bA_5$ singular points of~$B$, producing
the set of singularities
$$
\bD_4\splus\bA_5\splus2\bA_4.
$$
(Since $B$ splits into three components, we need to perturb two
points to get an irreducible curve.)
Assume that the points perturbed are those over~$R_5$. Then
the extra relations for $\pi_1(\Cp2\sminus B')$ are $\Ga=\Gb$ and
$\bGa=\bGb$, see Lemma~\ref{pert.a5},
and it is straightforward that the resulting group is abelian.

%\midinsert
%\plot{a11+a5+a3}
%\figure
%The set of singularities
%$(\bA_{11}\splus\bA_5)\splus\bA_3$
%\endfigure\label{fig.a11+a5+a3}
%\endinsert

\subsection{The set of singularities
$(\bA_{11}\splus\bA_5)\splus\bA_3$}\label{s.a11+a5+a3}
Take for~$\L$ be the section given by~\eqref{eq.a11+a5+a3}, see the
dotted grey line in Figure~\ref{fig.slanted}.
Choosing the generators
$(\eta_1,\eta_2,\eta_3,\eta_4)=(\Ga,\Gb,\Gd,\Gg)$
%$\Ga$, $\Gb$, $\Gd$, $\Gg$
in a real
generic fiber between~$R_5$ and~$R_1$ (\eg, over $x=1$), we obtain
the following relations:
$$
\alignat2
&[\Gb,\Gd\Gg]=[\Gd,\Gg\Gb]=[\Gg,\Gb\Gd]=1&\quad&
 \text{(the fiber through~$R_1$)},\\\allowdisplaybreak
&[\Gd,\Ga\Gb]=1&&
 \text{(the fiber through~$R_5$)},\\\allowdisplaybreak
&\Gd(\Ga\Gb)^3=\Gb\Ga\Gb\Gd\Ga\Gb\Ga&&
 \text{(the fiber through~$R_5$)},\\\allowdisplaybreak
&(\Gb\Ga\Gb\Gd)\1\Ga(\Gb\Ga\Gb\Gd)=\Gg&&
 \text{(the vertical tangent)},\\\allowdisplaybreak
&(\Ga\Gb\Gd\Gg)^2=1&&
 \text{(the relation at infinity)}.
\endalignat
$$
Passing to $\Ga$, $\bGa$, $\Gb$, $\bGb$, $\Gg$, $\bGg$, see
Lemma~\ref{bpi->pi}, we obtain a presentation for the group of a
reducible sextic of torus type with the set of singularities
$(\bA_{11}\splus\bA_5)\splus\bA_3$:
$$
\gather
\bGb\Gg=\Gg\Gb=\bGg\bGb=\Gb\bGg,
 \eqtag\label{eq.8.1}\\\allowdisplaybreak
\Ga\Gb=\bGa\bGb,
 \eqtag\label{eq.8.2}\\\allowdisplaybreak
(\Ga\Gb)^3=\bGb\bGa\bGb\Ga\Gb\Ga=\Gb\Ga\Gb\bGa\bGb\bGa,
 \eqtag\label{eq.8.3}\\\allowdisplaybreak
(\Gb\Ga)\1\Ga(\Gb\Ga)=\Gg,\quad
 (\bGb\bGa)\1\bGa(\bGb\bGa)=\bGg,
 \eqtag\label{eq.8.4}\\\allowdisplaybreak
\Gg\Ga\Gb\bGg\bGa\bGb=1.
 \eqtag\label{eq.8.5}
\endgather
$$

The perturbation $\bA_5\mapsto\bA_4$ produces an irreducible
plane sextic with the set of singularities
$$
\bA_{11}\splus\bA_4\splus\bA_3,
$$
while adding the relations
$\Ga=\Gg\1\bGg\bGa\bGg\1\Gg=\Gg=\Gg\1\bGg\Gg$,
see Lemma~\ref{pert.a5}, which imply $\Gg=\bGg=\Ga=\bGa$.
(Observe that the standard generators in
a fiber over $x\gg0$ are $\Ga_1=\Ga$, $\Gb_1=\Gg$,
$\Gd_1=\Gg\1\Gd\Gg$, $\Gg_1=\Gb$, and the extra relations given by
Lemma~\ref{pert.a5} are
$\Ga_1=\Gd_1\1\Ga_1\Gd_1=\Gb_1=\Gd_1\1\Gb_1\Gd_1$.)
The perturbations $\bA_{11}\mapsto\bA_{10}$
or $\bA_6\splus\bA_4$ produce irreducible sextics with the sets of
singularities
$$
\bA_{10}\splus\bA_5\splus\bA_3,\quad
 \bA_6\splus\bA_5\splus\bA_4\splus\bA_3
$$
and add the relations $\Ga=\bGa=\Gb=\bGb$, see
Lemma~\ref{pert.a11}. In both cases, it is immediate that
the resulting group is abelian.

\section{Digression: classification of reducible sextics}\label{S.reducible}

The curves considered in \S\ref{S.groups} are
sextics of torus type splitting into a quartic and a
conic. Here, we state and indicate the proofs of a few results
concerning the classification and the fundamental groups of such
curves. Details
%of the proofs
will be published elsewhere. In Section~\ref{proof.moduli.3}, we
outline the proof of Theorem~\ref{th.moduli.3}.

\subsection{The symmetry}
Theorem~\ref{splitting.symmetry} below substantiates
Conjecture~4.2.3 in~\cite{symmetric}, concerning the relation
between involutive stable symmetries of plane sextics and maximal
trigonal curves in~$\Sigma_2$.

\theorem\label{splitting.symmetry}
Let~$B$ be a plane sextic of torus type, with simple singularities
only, splitting into
\dashes
\dash
an irreducible quartic and irreducible conic,
\dash
an irreducible quartic and two lines, or
\dash
three irreducible conics.
\enddashes
Then $B$ admits an involutive stable \rom(in the sense
of~\cite{symmetric}\rom) symmetry $c\:\Cp2\to\Cp2$, and the
quotient $B/c$ is the maximal trigonal curve $\B_2\subset\Sigma_2$
with the set of singularities $\bA_5\splus\bA_2\splus\bA_1$, see
Section~\ref{s.a5+a2}.

Conversely, for any section~$\L$
of~$\Sigma_2$ not tangent to~$\B_2$ at its type~$\bA_5$ singular
point, the double $B=\DC(\B_2,\L)$ is a plane sextic
of torus type splitting as above.
\endtheorem

\Remark
The condition that $\L$ should not be tangent to~$\B_2$ at its
type~$\bA_5$ singular point is necessary and sufficient for~$B$ to
have simple singularities only.
\endRemark

\Remark\label{rem.S3}
For most sextics~$B$ as in Theorem~\ref{splitting.symmetry}, the
group of stable symmetries of~$B$ is~$\CG2$. Exceptions are
sextics splitting into three irreducible conics: for each such
sextic~$B$, the group of stable symmetries of~$B$ is the group~$\SG3$ of
permutations of the components of~$B$.
\endRemark

\proof
The proof is similar to~\cite{degt.Oka3}, \cite{degt.8a2},
and~\cite{symmetric}.
Assume that $B$ splits into a quartic~$B_4$ and a conic~$B_2$. It is clear
that the inner singularities are two cusps~$R'_\infty$, $R''_\infty$
of~$B_4$ (which may
degenerate to a single type~$\bA_5$ or~$\bE_6$ singular
point~$R_\infty$
of~$B_4$) and two type~$\bA_5$ points~$R'_5$, $R''_5$
of inflection tangency
of~$B_4$ and~$B_2$ (which may degenerate to a single
type~$\bA_{11}$ point~$R_5$ of $6$-fold intersection). The outer
singularities are the two points~$R_1'$, $R_1''$
of simple intersection of~$B_4$
and~$B_2$, which may degenerate to a single type~$\bA_3$
point~$R_1$.
Besides, $B_4$ may have an extra node or cusp~$Q_4$, and
$B_2$ may have an
extra node~$Q_2$ (splitting into two lines). As a further
degeneration, $R_1$ may merge with~$Q_4$ or~$Q_2$.
%and $B_2$ may pass
%through this point, replacing the two nodes above with a single
%type~$\bD_4$ or~$\bD_5$ singular point; $B_2$ may also have an
%extra node (splitting into two lines), and $B_4$ may pass through
%this node, forming a type~$\bD_4$ singular point.

Consider the minimal resolution~$\tX$ of the double covering
$X\to\Cp2$ ramified at~$B$. It is a $K3$-surface. Let
$L=H_2(\tX)$, let $\Sigma_P\subset L$ be the resolution
lattice of a singular point~$P$ of~$B$, and let~$\Gamma_P$ be the
Dynkin graph of~$P$.
Denote $\Sigma=\bigoplus_P\Sigma_P$ and
$\Gamma=\bigcup_P\Gamma_P$, and consider the lattice
$S=\Sigma\splus\<h>\subset L$, where $h\in L$ is the class of the
pull-back of a generic line in~$\Cp2$. One has $h^2=2$.
Let~$\CK$ be the image of $L=L^*$ in
the discriminant group $\discr S=S^*\!/S$.
Since $B$ is of torus
type, $\CK$ has an element of order~$3$, see~\cite{JAG}.
Besides, $\CK$ has an element of order~$2$ responsible for the
splitting $B=B_4+B_2$, see~\cite{degt.Oka}. (For example, the
element of order~$2$ is represented by the class
$[B_2]\in H_2(X)=L$.)

Consider an involutive symmetry $c_\Gamma\:\Gamma\to\Gamma$
acting as follows: $c_\Gamma$ transposes the two points within
each pair $(R'_\infty,R''_\infty)$, $(R'_5,R''_5)$, $(R'_1,R''_1)$
and acts identically on
the diagram of each point~$Q_\sdot$ that
does not coincide with~$R_1$. If a pair of points $R'_\sdot$,
$R''_\sdot$ merges to a single point~$R_\sdot$, then $c_\Gamma$
acts on $\Gamma_{\!R_\sdot}$ by its nontrivial symmetry.
(This description determines~$c_\Gamma$ uniquely up to a symmetry
of the diagrams of $R''_\infty$ and $R''_5$ whenever these points
are separate.)
Let $c_S\:S\to S$ be the extension of~$c_\Gamma$ identical on~$h$.
One can check that $c_\Gamma$ can be chosen so that
$c_S$ preserves~$\CK$ and induces the
identity on $\CK^\perp\!/\CK$; hence, $c_S$ extends to an
involutive automorphism $c_*\:L\to L$ identical on~$\Sigma^\perp$.
The latter is induced by a unique involution $c\:\Cp2\to\Cp2$
preserving~$B$, \cf.~\cite{symmetric}. Details are left to the
reader.

The fact that the quotient $B/c$ is the curve~$\B_2$ is
straightforward: the quotient must have singular points of
types~$\bA_5$, $\bA_2$, and~$\bA_1$, resulting from the (pairs of)
points~$R_5$, $R_\infty$, and~$R_1$, respectively, and, due
to~\cite{symmetric}, such a curve is unique. The converse
statement is obvious.
\endproof

If $B$ splits into three conics, it has three type~$\bA_5$ inner
singular points and three outer nodes, which may merge to a single
type~$\bD_4$ singular point. An order~$3$ symmetry of~$B$, see
Remark~\ref{rem.S3}, is constructed as above, starting from the
order~$3$ symmetry $c_\Gamma\:\Gamma\to\Gamma$ transposing
cyclically the three inner points and three nodes (or acting by an
order three symmetry on the Dynkin graph~$\bD_4$).

\subsection{The classification}
Using the stable symmetry and the description of special sections
found in~\ref{s.a5+a2} and~\ref{s.others}, one immediately
obtains a deformation classification of sextics splitting as in
Theorem~\ref{splitting.symmetry}.

\theorem\label{splitting.moduli}
Let~$B$ be a sextic as in Theorem~\ref{splitting.symmetry}.
Then the combinatorial type of singularities of~$B$ is one of
those listed in~\cite{Oka.reducible} or
$(\bE_6\splus\bA_{11})\splus2\bA_1$.
The equisingular moduli space of sextics as in
Theorem~\ref{splitting.symmetry} realizing each of these
combinatorial types is rational\rom; in particular, it is
connected.
\endtheorem

\Remark
When referring to~\cite{Oka.reducible}, we mean the combinatorial
types marked as $B_2+B_4$, $B_1+B_1'+B_4$, or $B_2+B_2'+B_2''$ in
Theorem~$1$.
The set of singularities
$(\bE_6\splus\bA_{11})\splus2\bA_1$ (an irreducible quartic
with a type~$\bE_6$ singular point and a conic,
see Section~\ref{s.a11+e6})
is erroneously missing in~\cite{Oka.reducible}.
\endRemark

\proof
As in Section~\ref{proof.moduli}, the equisingular moduli spaces
$\CM(\Sigma)$ of sextics splitting as in
Theorem~\ref{splitting.symmetry} and
possessing a given set of singularities~$\Sigma$
(more precisely, the spaces $\tilde\CM(\Sigma)$ of pairs $(B,c)$,
where $B$ is a sextic and $c$ is an involutive stable symmetry
of~$B$) can be identified with the spaces of sections~$\L$
of~$\Sigma_2$ that are in a certain prescribed special position
with respect to~$\B_2$. The latter are described in
Sections~\ref{s.a5+a2} and~\ref{s.others}; they are all rational.
It remains to notice that the forgetful map
$\tilde\CM(\Sigma)\to\CM(\Sigma)$ is one to one, as any two stable
involutions of a sextic~$B$ are projectively equivalent, see
Remark~\ref{rem.S3}.
\endproof

\theorem\label{splitting.group}
Let~$B$ be a sextic as in Theorem~\ref{splitting.symmetry}. Then
the fundamental group $\pi_1(\Cp2\sminus B)$ factors to the group
given by~\eqref{eq.D4.perturbed}.
\endtheorem

\Remark\label{rem.group}
The groups of most maximal sextics~$B$ as in
Theorem~\ref{splitting.symmetry} are computed in
Sections~\ref{s.a11+3a2}--\ref{s.a11+a5+a3}. Perturbing~$\L$, one
can easily see that, if $B$ has exactly two components, the
quartic component of~$B$ has a set of singularities other
than~$3\bA_2$, and $\mu(B)<19$, then $\pi_1(\Cp2\sminus B)$ is the
group given by~\eqref{eq.D4.perturbed}.
\endRemark

\proof
Any sextic as in Theorem~\ref{splitting.symmetry} can be perturbed
to a `simplest' one, with the set of singularities
$(2\bA_5\splus2\bA_2)\splus2\bA_1$, which is the double of~$\B_2$
ramified at a section~$\L$ transversal to~$\B_2$.
The group of a simplest sextic is~\eqref{eq.D4.perturbed}, see,
\eg, Section~\ref{s.a11+2a2+d4}.
\endproof

\subsection{Proof of Theorem~\ref{th.moduli.3}}\label{proof.moduli.3}
Let~$P$ be the type~$\bA_{17}$ singular point, and
let $\Sigma_P\subset\Sigma\subset S\subset L$ \etc. be as in the
proof of Theorem~\ref{splitting.symmetry}. Denote
$S_P=\Sigma_P\oplus\<h>$.
Since the sextic is
reducible and of torus type, the intersection
$\CK'=\CK\cap\discr S_P$ contains an element of
order~$2$ and an element of order~$3$, see~\cite{degt.Oka}
and~\cite{JAG}. On the other hand, $\ls|\discr S_P|=36$; hence,
$(\CK')^\perp\!/\CK'=0$, \ie, the primitive hull of~$S_P$ in~$L$
is unimodular and the classification of homological types
reduces to the study of sublattices isomorphic to~$0$, $\bA_1$,
$2\bA_1$, or~$\bA_2$ in the direct sum of two hyperbolic planes.
The rest is straightforward.
\qed

\Remark
From the proof, it follows that each sextic~$B$ as in
Theorem~\ref{th.moduli.3} admits a stable involutive symmetry~$c$
(constructed as in the proof of Theorem~\ref{splitting.symmetry}
starting from the nontrivial symmetry of~$\Gamma_P$). However, one
has $O_c\in B$; hence, the quotient $B/c\subset\Sigma_2$ is not a
trigonal curve but rather a hyperelliptic curve with a
type~$\bA_7$ singular point at~$E$. Thus, Conjecture~4.2.3
in~\cite{symmetric} needs to be modified to include maximal, in
some sense, hyperelliptic curves as well.
\endRemark

\section{Other fundamental groups}\label{S.others}

We consider a triple $\B_2'$, $\L'$, $\L$ as in~\S\ref{S.groups}
and make~$\L$ and~$\L'$ trade r\^oles, \ie, we let $\B=\B_2'+\L$ and
consider the double $B=\DC((\B_2'+\L),\L')$ ramified at~$\L'$.
The groups
$\pi_1(\Sigma_2\sminus(\B_2'\cup\L'\cup\L\cup E))$ are computed
in~\S\ref{S.groups}, and we merely use an appropriately modified
version of Lemma~\ref{bpi->pi} (with the r\^ole of~$\Gd$ played by
the generator corresponding to~$\L'$) to obtain
$\pi_1(\Cp2\sminus B)$. Then, as in~\S\ref{S.groups}, we
perturb~$B$ to an irreducible sextic~$B'$ and apply
Lemmas~\ref{pert.a11}--\ref{pert.e7}.

\subsection{The set of singularities
$2\bE_7\splus\bD_5$}\label{s.2e7+d5}
Take for~$\L$ the section given by~\eqref{eq.a11+e6}, see
Section~\ref{s.a11+e6} and Figure~\ref{fig.horizontal} (the solid
grey line).
The resulting sextic~$B$ has
the set of singularities $2\bE_7\splus\bD_5$, and
the group $\pi_1(\Cp2\sminus B)$ is obtained from~\eqref{eq.2E7.1}
by letting $\Gb^2=1$ and passing to the subgroup generated by $\Ga$,
$\bGa=\Gb\Ga\Gb$, $\Gd$, and $\bGd=\Gb\Gd\Gb$. One has $\Gd=\bGd$
and hence
$$
\pi_1(\Cp2\sminus B)=\bigl<\Ga,\bGa,\Gd\big|
 \Ga\bGa\Ga=\bGa\Ga\bGa,\
 [\Ga,\Gd]=[\bGa,\Gd]=1,\
 \Ga^2\bGa^2\Gd^2=1\bigr>.
\eqtag\label{eq.minimal}
$$

The irreducible perturbation $\bD_5\mapsto\bA_4$ produces the set
of singularities
$$
2\bE_7\splus\bA_4
$$
and adds the relation $\Ga=\bGa=\Gd$, see Lemma~\ref{pert.d5}. The
irreducible
perturbations $\bE_7\mapsto\bE_6$, $\bA_6$, or $\bA_4\splus\bA_2$
of one of the two type~$\bE_7$ singular points of~$B$
produce the sets of
singularities
$$
\bE_7\splus\bE_6\splus\bD_5,\quad
 \bE_7\splus\bD_5\splus\bA_6,\quad
 \bE_7\splus\bD_5\splus\bA_4\splus\bA_2
$$
while adding at least the relation $\Ga\Gd\Ga=\Gd\Ga\Gd$ (or
$\Gg\Gd\Gg=\Gd\Gg\Gd$, where $\Gg=\Ga\1\bGa\Ga$), see
Lemma~\ref{pert.e7} and~\eqref{eq.e7}.
In each case, it is immediate that the
resulting fundamental group $\pi_1(\Cp2\sminus B')$ is abelian.

\subsection{The set of singularities
$2\bE_7\splus\bA_3\splus\bA_2$}\label{s.2e7+a3+a2}
Take for~$\L$ the section given by~\eqref{eq.2a5+e6+a3}, see
%Section~\ref{s.2a5+e6+a3} and
%Figure~\ref{fig.horizontal} (the dotted grey line).
the dotted grey line in Figure~\ref{fig.horizontal}.
A presentation for the group
$\pi_1(\Sigma_2\sminus(\B\cup\L'\cup E))$ is found in
Section~\ref{s.2a5+e6+a3}. To pass to $\pi_1(\Cp2\sminus B)$, we
let $\Gb^2=1$ and consider the subgroup generated by~$\Ga$,
$\bGa=\Gb\Ga\Gb$, $\Gg$, $\bGg=\Gb\Gg\Gb$,
$\Gd$, and $\bGd=\Gb\Gd\Gb$. The relations are
$$
\gather
\Ga\bGa\Ga=\bGa\Ga\bGa,
 \eqtag\label{eq.9.1}\\\allowdisplaybreak
\Ga\bGa\Ga\1=\Gg,\quad
 \bGa\Ga\bGa\1=\bGg,
 \eqtag\label{eq.9.2}\\\allowdisplaybreak
\Ga\1\Gd\Ga=\bGa\1\bGd\bGa,
 \eqtag\label{eq.9.3}\\\allowdisplaybreak
\Gd\bGg=\Gg\Gd=\bGg\bGd=\bGd\Gg,
 \eqtag\label{eq.9.4}\\\allowdisplaybreak
\Gg\Gd\Ga\bGg\bGd\bGa=1.
 \eqtag\label{eq.9.5}
\endgather
$$

The irreducible perturbation $\bA_3\mapsto\bA_2$ produces the set
of singularities
$$
2\bE_7\splus2\bA_2
$$
and adds the relation $\Gg=\bGg=\Gd=\bGd$. Then,
comparing~\eqref{eq.9.2} and~\eqref{eq.9.3}, one concludes that
$\Ga=\bGa$ and hence $\Ga=\Gg$. Thus, the group is abelian.

Consider a maximal irreducible perturbation of one of the two
type~$\bE_7$ singular points of~$B$,
producing irreducible sextics with the sets of singularities
$$
\bE_7\splus\bE_6\splus\bA_3\splus\bA_2,\quad
 \bE_7\splus\bA_6\splus\bA_3\splus\bA_2,\quad
 \bE_7\splus\bA_4\splus\bA_3\splus2\bA_2,
$$
see Lemma~\ref{pert.e7}.
A generic real fiber close
to~$R_\infty$ (the
type~$\bE_7$ singular point of~$\B$) is over $x\gg0$, and the
standard generators in this fiber
are $\Ga$, $(\Gb\Gg\Gd)\Gd(\Gb\Gg\Gd)\1=\Gd$, $\Gb\Gg\Gb\1$,
and~$\Gb$. Hence, the group $\pi_1(\MB\sminus B)$ of a Milnor ball
about the point perturbed
is generated by $\Ga$, $\Gd$, and~$\bGg$, and the perturbation
adds at least the relation $\bGg\Ga\bGg=\Gd\bGg\Ga$ (the third
relation in~\eqref{eq.e7}\,).
Using~\eqref{eq.9.4}, the additional relation simplifies to
$\Ga\bGg=\bGd\Ga$. Hence, $\bGd=\Ga\bGg\Ga\1=\bGa$
%(from~\eqref{eq.9.2} using~\eqref{eq.9.1}),
(substituting to~\eqref{eq.9.2} and using~\eqref{eq.9.1}\,),
$\Ga\1\Gd\Ga=\bGa$ (from~\eqref{eq.9.3}\,),
$\Gd=\Ga\bGa\Ga\1=\Gg$ (from~\eqref{eq.9.2}\,),
$\Gg=\bGg$ and $[\bGg,\bGa]=1$ (from~\eqref{eq.9.4} again),
and $\Ga=\bGg=\Gg=\bGa$ (from~\eqref{eq.9.2}\,).
Thus, the group is abelian.

\subsection{The set of singularities
$2\bE_7\splus\bA_2\splus3\bA_1$}\label{s.2e7+a2+3a1}
Consider the section~$\L$ tangent to~$\L'$ and tangent to~$\B_2'$
at its cusp~$R_\infty$. It is given by $y=3/4$.
Choose the generators
%$\Ga$, $\Gd$, $\Gb$, $\Gg$
$(\eta_1,\eta_2,\eta_3,\eta_4)=(\Ga,\Gd,\Gb,\Gg)$
in a real fiber~$F$ over~$x\gg0$.

We are only interested in the
perturbation $\bE_7\mapsto\bA_6$ of \emph{both} type~$\bE_7$
singular points of~$B$, producing an irreducible sextic with the
set of singularities
$$
2\bA_6\splus\bA_2\splus3\bA_1.
$$
This perturbation can be realized symmetrically, by perturbing the
type~$\bE_7$ singular point of~$\B$ in~$\Sigma_2$. According to
Lemma~\ref{pert.e7}, this gives the relations $\Ga=\Gd=\Gb$, and
the monodromy about~$R_1$ adds the
relation $[\Gb,\Gg]=1$. Hence, the resulting group
$\pi_1(\Cp2\sminus B')$ is abelian.

\subsection{The set of singularities
$2\bD_5\splus\bA_7\splus\bA_2$}\label{s.2d5+a7+a2}
Let~$\L$ be the section given by~\eqref{eq.3a5+d4},
see the solid grey line in Figure~\ref{fig.slanted}.
As in Section~\ref{s.3a5+d4}, the
relations for $\pi_1(\Cp2\sminus B)$ are
\eqref{eq.9.1}--\eqref{eq.9.3}, \eqref{eq.9.5}, and the relations
$$
\Gg\Gd=\bGg\bGd,\quad
 \bGd\Gg\Gd\Gg=(\Gg\Gd)^2,\quad
 \Gd\bGg\bGd\bGg=(\bGg\bGd)^2
 \eqtag\label{eq.10.1}
$$
replacing~\eqref{eq.9.4}.

The irreducible perturbations $\bA_7\mapsto\bA_6$ or
$\bA_4\splus\bA_2$, see Lemma~\ref{pert.a7},
produce the sets of singularities
$$
2\bD_5\splus\bA_6\splus\bA_2,\quad
 2\bD_5\splus\bA_4\splus2\bA_2
$$
and add the relations $\Gg=\bGg=\Gd=\bGd$. As in
Section~\ref{s.3a5+d4}, the resulting groups are abelian.
The irreducible perturbation $\bD_5\mapsto\bA_4$ of one of the
type~$\bD_5$ singular points produces the set of singularities
$$
\bD_5\splus\bA_7\splus\bA_4\splus\bA_2.
$$
The standard generators in a generic fiber close to~$R_\infty$
(a fiber over~$x\gg0$) are $\Ga$,
$(\Gb\Gg\Gd)\Gg(\Gb\Gg\Gd)\1$, $(\Gb\Gg)\Gd(\Gb\Gg)\1$, and~$\Gb$.
In view of Lemma~\ref{pert.d5}, the extra relations added to
the group are $\Ga=(\bGg\bGd)\bGg(\bGg\bGd)\1=\bGg\bGd\bGg\1$.
This implies $\Ga=\bGg=\bGd$, and using
\eqref{eq.9.1}--\eqref{eq.9.3} one derives that
$\Ga=\bGa=\Gg=\Gd$. Hence, the group is abelian.

\subsection{The set of singularities
$3\bD_5\splus\bA_3$}\label{s.3d5+a3}
Let~$\L$ be the section given by~\eqref{eq.a11+a5+a3},
see the dotted grey line in Figure~\ref{fig.slanted}.
Starting from the presentation found in Section~\ref{s.a11+a5+a3},
letting $\Gb^2=1$, and passing to the subgroup generated by $\Ga$,
$\bGa=\Gb\Ga\Gb$, $\Gg$, $\bGg=\Gb\Gg\Gb$, $\Gd$,
$\bGd=\Gb\Gd\Gb$, we obtain the following relations for
%the group
$\pi_1(\Cp2\sminus B)$:
$$
\gather
\Gd\Gg=\bGd\bGg=\Gg\bGd=\bGg\Gd,
 \eqtag\label{eq.11.1}\\\allowdisplaybreak
\Gd\Ga=\Ga\bGd,\quad
 \bGd\bGa=\bGa\Gd,
 \eqtag\label{eq.11.2}\\\allowdisplaybreak
\Gd\Ga\bGa\Ga=\bGd\bGa\Ga\bGa,
 \eqtag\label{eq.11.3}\\\allowdisplaybreak
\Ga\1\bGa\Ga=\Gg,\quad
 \bGa\1\Ga\bGa=\bGg,
 \eqtag\label{eq.11.4}\\\allowdisplaybreak
\Gd\Gg\Ga\bGd\bGg\bGa=1.
 \eqtag\label{eq.11.5}
\endgather
$$
The irreducible perturbation $\bA_3\mapsto\bA_2$ produces the set
of singularities
$$
3\bD_5\splus\bA_2,
$$
adding the relation $\Gg=\bGg=\Gd=\bGd$.
%, see Lemma~\ref{pert.a5}.
The irreducible perturbation $\bD_5\mapsto\bA_4$ of one of the
type~$\bD_5$ singular points of~$B$ produces the set of
singularities
$$
2\bD_5\splus\bA_4\splus\bA_3,
$$
adding the relation $\Ga=\bGa=\Gd=\bGd$, see Lemma~\ref{pert.d5}.
(We can assume that the point perturbed is over~$R_5$.) It is
straightforward that, in both cases, the new fundamental
group is abelian.

\subsection{Concluding remarks}\label{s.remarks}
The groups of all reducible curves obtained in this section are
non-abelian; they all factor to the `minimal' group~$G$ given
by~\eqref{eq.minimal}, which can also be described as a central
extension
$$
1@>>>\<\Gd>@>>>G@>>>\AG4@>>>1.
$$
This result is quite expectable, as all curves split into a
conic~$B_2$ and a three cuspidal quartic~$B_4$, and
$\pi_1(\Cp2\sminus B_4)=\AG4$.

It is worth mentioning that all curves admit regular
$\SG3$-coverings while obviously not being of torus type. Hence,
Theorem~4.1.1 in~\cite{degt.Oka} does not extend to reducible
curves literally. Certainly, the reason is the fact that the
cyclic part of the covering is not ramified at~$B$ but rather
at~$B_4$ only.

In most cases, the trigonal curve $\B=\B_2'+\L\subset\Sigma_2$
used in the construction is maximal, with the set of
singularities $\bE_7\splus\bA_1$ or $\bD_5\splus\bA_3$. Thus, one
may hope that the deck translation is a stable symmetry of the
covering sextic~$B$ (\cf. Conjecture~4.2.3 in~\cite{symmetric})
and use this correspondence to classify sextics.

In the former case, the set of singularities $\bE_7\splus\bA_1$,
the sextics are characterized by the splitting $B=B_4+B_2$,
where $B_4$ is a quartic with at least two cusps and $B_2$ is a
conic (possibly reducible) tangent to~$B_4$ at two of its cusps.
Any such curve is indeed symmetric: in appropriate affine
coordinates $(x,y)$ in~$\Cp2$, the curves~$B_4$ and~$B_2$ can be
given by $a+b(x+y)+cx^2y^2=0$ and $d+exy=0$, respectively.
There are $13$ deformation families of such curves, of which four
are maximal. Three maximal families are considered in
Sections~\ref{s.2e7+d5}--\ref{s.2e7+a2+3a1}; the fourth one has
the set of singularities $2\bE_7\splus\bD_4\splus\bA_1$ (the
conic~$B_2$ splits and the quartic~$B_4$ has an extra node and
passes through the node of~$B_2$).

In the latter case, the set of singularities $\bD_5\splus\bA_3$,
the sextic splits into $B_4+B_2'$, where $B_4$ is a quartic with at
least two cusps and $B_2'$ is a conic passing through two cusps
of~$B_4$ and tangent to~$B_4$ at all other intersection points. It
appears that this configuration, as well as most of its
degenerations, is realized by several equisingular deformation
families, only one of them admitting a stable symmetry. The
symmetric sextics seem to be related to the sextics of torus
type considered in~\S\ref{S.reducible}: they are obtained by
replacing the conic~$B_2$ in the splitting $B=B_4+B_2$, see
Theorem~\ref{splitting.symmetry}, by the conic $B_2'=\{p=0\}$,
where $p^3+q^2=0$ is the (only) torus structure on~$B$.
(From the point of view of the trigonal curve, we replace the
$\L'$ component in $\B_2=\B_2'+\L'$ with the section passing
through~$R_\infty$ and tangent to~$\B_2'$ at~$R_5$.)
We will
treat this subject in details elsewhere.

Note that each double $B=\DC(\B,\L)$ of the trigonal curve~$\B$
with the set of singularities $\bD_5\splus\bA_3$ has non-abelian
fundamental group: all groups factor to the `simplest' one,
corresponding to the case when $\L$ is transversal to~$\B$.
Letting $\Ga=\bGa$, $\Gg=\bGg$, and $\Gd=\bGd$ in the presentation
in Section~\ref{s.2d5+a7+a2}, we obtain the following presentation
for the simplest group~$G$:
$$
G=\<\Gg,\Gd\,|\,(\Gg\Gd)^2=(\Gd\Gg)^2,\ (\Gg\Gd\Gg)^2=1>.
$$
Introducing new generators $u=\Gd\Gg$, $v=\Gg\Gd\Gg$, we can
simplify the presentation to $G=\<u,v\,|\,v^2=[v,u^2]=1>$. It is
clear that the commutant $[G,G]$ equals~$\Z$, both~$u$ and~$v$
acting on $[G,G]$ by the multiplication by~$(-1)$. In particular,
all curves admit regular $\DG{2n}$-coverings for any $n\ge3$;
however, they are \emph{not} \term2n-sextics.

\section{Summary}\label{S.summary}

We summarize the results on the fundamental group of an
irreducible plane sextic obtained in this
paper and combine them with~\cite{degt.8a2} and~\cite{degt.e6}. We
confine ourselves to the case of simple singularities only; the
groups of sextics with a non-simple singular point are essentially
found in~\cite{degt.Oka} and~\cite{degt.Oka2} (see
also Oka, Pho~\cite{OkaPho}).

\midinsert
\table
Sextics of torus type: known groups
\endtable\label{tab.known}
\centerline{\vbox{\offinterlineskip\def\ex{\omit\vrule height2pt&&\cr}
\def\*{\hbox to0pt{\hss$^*$}}
\halign{\tabstrut\quad\hss#\hss\quad&\vrule
 \quad$#$\hss\quad\vrule&\quad#\hss\quad\vrule\cr
\noalign{\hrule}\ex
No.&\ \text{The set of singularities}&\ Where?\cr\ex
\noalign{\hrule}\ex
nt23&\*(6\bA_2)\splus3\bA_2&see \cite{degt.8a2}\cr
nt32&\*(\bA_5\splus4\bA_2)\splus\bE_6&see \cite{degt.8a2}\cr
nt36&(\bA_5\splus4\bA_2)\splus\bA_4\splus\bA_1&$\to$ nt32\cr
nt47&\*(\bE_6\splus4\bA_2)\splus\bA_5&same as nt32\cr
nt54&(2\bA_5\splus2\bA_2)\splus\bA_1&$\to$ nt55\cr
nt57&(2\bA_5\splus2\bA_2)\splus\bA_4&see \ref{s.2a5+2a2+d5}\cr
nt63&(\bE_6\splus\bA_5\splus2\bA_2)\splus\bA_3&$\to$ nt70\cr
nt67&\*(\bE_6\splus\bA_5\splus2\bA_2)\splus2\bA_2&same as nt32\cr
nt70&\*(2\bE_6\splus2\bA_2)\splus\bA_3&see \cite{degt.e6}\cr
nt74&(\bA_8\splus3\bA_2)\splus\bA_3&$\to$ nt128\cr
nt77&(\bA_8\splus3\bA_2)\splus\bD_4&see \ref{s.a11+2a2+d4}\cr
nt81&(\bA_8\splus3\bA_2)\splus\bA_2\splus\bA_1&$\to$ nt88\cr
nt88&(\bA_8\splus3\bA_2)\splus\bA_2\splus2\bA_1&see \ref{s.a11+3a2}\cr
nt99&\*(2\bE_6\splus\bA_5)\splus\bA_2&see \cite{degt.e6}, \cite{EyralOka}\cr
nt100&\*(3\bE_6)\splus\bA_1&see \cite{degt.e6}, \cite{OkaPho}\cr
nt103&(\bA_8\splus\bA_5\splus\bA_2)\splus\bA_3&$\to$ nt128\cr
nt105&(\bA_8\splus\bA_5\splus\bA_2)\splus2\bA_1&$\to$ nt103\cr
nt108&(\bE_6\splus\bA_8\splus\bA_2)\splus\bA_1&$\to$ nt112\cr
nt112&(\bE_6\splus\bA_8\splus\bA_2)\splus2\bA_1&see \ref{s.a11+e6}\cr
nt117&(\bA_{11}\splus2\bA_2)\splus\bA_3&see \ref{s.a11+2a2+d4}\cr
nt120&(\bA_{11}\splus2\bA_2)\splus\bA_2\splus\bA_1&see \ref{s.a11+3a2}\cr
nt128&(2\bA_8)\splus\bA_3&see \ref{s.2a8+a3}\cr
nt134&(\bA_{11}\splus\bA_5)\splus\bA_2&$\to$ nt145\cr
nt135&(\bE_6\splus\bA_{11})\splus\bA_1&see \ref{s.a11+e6}\cr
nt138&(\bA_{14}\splus\bA_2)\splus\bA_2&$\to$ nt145\cr
nt145&(\bA_{17})\splus\bA_2&see \ref{s.a17+a2}\cr
\ex\noalign{\hrule}
\crcr}}}
\endinsert

\subsection{Sextic of torus type}\label{s.torus}
According to Oka, Pho~\cite{OkaPho.moduli},
there are $19$ tame and $109$ non-tame sets of simple
singularities realized by irreducible plane sextics of torus type.
At present, the fundamental group is known for $113$
sets of singularities, including all tame ones. The result is
summarized in Table~\ref{tab.known}, where `nt\#' is the
notation introduced in~\cite{OkaPho.moduli} and the last column
indicates a proof, either by referring to the appropriate
paper/section or by suggesting a degeneration (in the form
`\hbox{$\to$ nt\#}')
to a set of singularities with known group. We only list the sets
of singularities for which the degenerations suggested
in~\cite{OkaPho.moduli} lead to sextics whose groups are still
unknown.

With few exceptions, the fundamental group of an irreducible
sextic of torus type is Zariski's group $\RBG3\cong\CG2*\CG3$. The
known exceptions are
\dashes
\dash
sextics of weight~$8$ and~$9$ in the sense of~\cite{degt.Oka},
see~\cite{degt.8a2};
\dash
sextics marked with a * in Table~\ref{tab.known}, see references
in the table;
\dash
the set of singularities $2\bE_6\splus2\bA_2\splus2\bA_1$,
see~\cite{degt.e6}.
\enddashes
Various perturbations of the exceptional sextics are studied
explicitly in~\cite{degt.8a2} and~\cite{degt.e6}; all other groups
are given by Corollary~\ref{torus->torus}.

\Remark
For most non-maximal sets of singularities, the connectedness of
the equisingular deformation family is still unknown, although
expected, see Conjecture~\ref{conjecture}.
For these sets of singularities, we can only state the
result in the form of existence, \ie, to assert that there is a
sextic~$B$ of torus type realizing a given set of singularities
and such that $\pi_1(\Cp2\sminus B)=\RBG3$. To my knowledge, the
sets of singularities for which the classification is completed
are:
\dashes
\dash
sextics admitting a stable involutive symmetry, see~\cite{symmetric}
for the list and \cite{degt.8a2}, \cite{degt.e6}, and
Theorem~\ref{th.moduli} for the classification;
\dash
the sets of singularities of the form
$(\text{inner points})\splus k\bA_1$, see~\cite{JAG}
and~\cite{Aysegul}.
\enddashes
\endRemark

The fifteen remaining sets of singularities,
for which the fundamental group is still unknown,
are listed in Table~\ref{tab.unknown} (with a reference to the
notation of~\cite{OkaPho.moduli}).

\midinsert
\table
Sextics of torus type: unknown groups
\endtable\label{tab.unknown}
\def\ex{\omit\vrule height2pt&\cr}
\def\header{\vtop\bgroup\halign\bgroup
 \tabstrut\ \hss##\hss\ &
 \vrule\ \ ##\hss\ \vrule\cr
 \noalign{\hrule}\ex
 No.&The set of singularities\cr\ex
 \noalign{\hrule}\ex}
\def\footer{\ex\noalign{\hrule}\crcr\egroup\egroup}
\offinterlineskip\centerline{\header
 nt64&$(\bE_6\splus\bA_5\splus2\bA_2)\splus\bA_4$\cr
 nt75&$(\bA_8\splus3\bA_2)\splus\bA_4$\cr
 nt78&$(\bA_8\splus3\bA_2)\splus\bD_5$\cr
 nt82&$(\bA_8\splus3\bA_2)\splus\bA_3\splus\bA_1$\cr
 nt83&$(\bA_8\splus3\bA_2)\splus\bA_4\splus\bA_1$\cr
nt104&$(\bA_8\splus\bA_5\splus\bA_2)\splus\bA_4$\cr
nt106&$(\bA_8\splus\bA_5\splus\bA_2)\splus\bA_2\splus\bA_1$\cr
nt109&$(\bE_6\splus\bA_8\splus\bA_2)\splus\bA_2$\cr\footer\ \header
nt110&$(\bE_6\splus\bA_8\splus\bA_2)\splus\bA_3$\cr
nt113&$(\bE_6\splus\bA_8\splus\bA_2)\splus\bA_2\splus\bA_1$\cr
nt118&$(\bA_{11}\splus2\bA_2)\splus\bA_4$\cr
nt136&$(\bE_6\splus\bA_{11})\splus\bA_2$\cr
nt139&$(\bA_{14}\splus\bA_2)\splus\bA_3$\cr
nt141&$(\bA_{14}\splus\bA_2)\splus2\bA_1$\cr
nt142&$(\bA_{14}\splus\bA_2)\splus\bA_2\splus\bA_1$\cr\footer}
\endinsert

%What is left:
%
%$(\bE_6\splus\bA_{11})\splus\bA_2$
%
%$(\bE_6\splus\bA_8\splus\bA_2)\splus\bA_2\splus\bA_1$ $\to$
%$(\bE_6\splus\bA_8\splus\bA_2)\splus\bA_2$ $\leftarrow$
%$(\bE_6\splus\bA_8\splus\bA_2)\splus\bA_3$
%
%$(\bE_6\splus\bA_5\splus2\bA_2)\splus\bA_4$
%
%$(\bA_{14}\splus\bA_2)\splus\bA_2\splus\bA_1$ $\to$
%$(\bA_{14}\splus\bA_2)\splus2\bA_1$ $\leftarrow$
%$(\bA_{14}\splus\bA_2)\splus\bA_3$
%
%$(\bA_{11}\splus2\bA_2)\splus\bA_4$
%
%$(\bA_8\splus\bA_5\splus\bA_2)\splus\bA_4$ $\to$
%$(\bA_8\splus\bA_5\splus\bA_2)\splus\bA_2\splus\bA_1$
%
%$(\bA_8\splus3\bA_2)\splus\bD_5$ $\to$
%$(\bA_8\splus3\bA_2)\splus\bA_3\splus\bA_1$
%
%$(\bA_8\splus3\bA_2)\splus\bA_4\splus\bA_1$ $\to$
%$(\bA_8\splus3\bA_2)\splus\bA_4$

\subsection{Sextics with abelian fundamental groups}\label{s.abelian}
In Table~\ref{tab.abelian}, we list the
sets of singularities realized by irreducible plane sextics
with abelian fundamental group, together with the references to the
sections where these curves are constructed. Combining these
results with~\cite{degt.8a2} and~\cite{degt.e6} and considering
all perturbations, we obtain $768$ sets of singularities not
covered by Nori's theorem~\cite{Nori}.

\midinsert
\table
Extremal sextics with abelian fundamental groups
\endtable\label{tab.abelian}
\def\ex{\omit\vrule height2pt&\omit\vrule\hss\vrule\cr}%
\def\header{\vtop\bgroup\halign\bgroup
 \tabstrut\quad##\hss\quad&
 \vrule\quad\hss\null\getref##\hss\quad\vrule\cr
 \noalign{\hrule}\ex
 The set of singularities&\omit\vrule\hss\ Where?\ \hss\vrule\cr\ex
 \noalign{\hrule}\ex}
\def\footer{\ex\noalign{\hrule}\crcr\egroup\egroup}
{\offinterlineskip\obeylines\let
\getline\hbox to\hsize{\hss\header
[2 e[7] + a[4], "s.2e7+d5"]\cr
[2 e[7] + 2 a[2], "s.2e7+a3+a2"]\cr
[e[7] + e[6] + d[5], "s.2e7+d5"]\cr
[e[7] + e[6] + a[3] + a[2], "s.2e7+a3+a2"]\cr
[e[7] + d[5] + a[6], "s.2e7+d5"]\cr
[e[7] + d[5] + a[4] + a[2], "s.2e7+d5"]\cr
[e[7] + a[6] + a[3] + a[2], "s.2e7+a3+a2"]\cr
[e[7] + a[4] + a[3] + 2 a[2], "s.2e7+a3+a2"]\cr
%[e[6] + a[11] + a[1], "s.a11+e6"]\cr
[e[6] + a[10] + 2 a[1], "s.a11+e6"]\cr
%[e[6] + a[8] + a[2] + 2 a[1], "s.a11+e6"]\cr
[e[6] + a[6] + a[4] + 2 a[1], "s.a11+e6"]\cr
[e[6] + a[5] + a[4] + a[3], "s.2a5+e6+a3"]\cr
[3 d[5] + a[2], "s.3d5+a3"]\cr
[2 d[5] + a[6] + a[2], "s.2d5+a7+a2"]\cr
[2 d[5] + a[4] + a[3], "s.3d5+a3"]\cr
[2 d[5] + a[4] + 2 a[2], "s.2d5+a7+a2"]\cr
[d[5] + a[7] + a[4] + a[2], "s.2d5+a7+a2"]\cr
[d[5] + a[5] + a[4] + 2 a[2], "s.2a5+2a2+d5"]\cr
[d[4] + a[10] + 2 a[2], "s.a11+2a2+d4"]\cr\footer\ \header
%[d[4] + a[8] + 3 a[2], "s.a11+2a2+d4"]\cr
[d[4] + a[6] + a[4] + 2 a[2], "s.a11+2a2+d4"]\cr
[d[4] + a[5] + 2 a[4], "s.3a5+d4"]\cr
[a[16] + a[2], "s.torus.pert"]\cr
[a[15] + a[2] + a[1], "s.torus.pert"]\cr
[a[13] + a[3] + a[2], "s.torus.pert"]\cr
[a[12] + a[4] + a[2], "s.torus.pert"]\cr
[a[11] + a[4] + a[3], "s.a11+a5+a3"]\cr
[a[11] + a[3] + 2 a[2], "s.a11+2a2+d4"]\cr
%[a[11] + 3 a[2] + a[1], "s.a11+3a2"]\cr
[a[10] + a[6] + a[2], "s.torus.pert"]\cr
[a[10] + a[5] + a[3], "s.a11+a5+a3"]\cr
[a[10] + 3 a[2] + 2 a[1], "s.a11+3a2"]\cr
[a[9] + a[7] + a[2], "s.torus.pert"]\cr
[a[8] + a[7] + a[3], "s.torus.pert"]\cr
[a[8] + a[6] + a[3] + a[1], "s.torus.pert"]\cr
[a[8] + a[4] + 2 a[3], "s.torus.pert"]\cr
%[a[8] + 4 a[2] + 2 a[1], "s.a11+3a2"]\cr
[2 a[6] + a[2] + 3 a[1], "s.2e7+a2+3a1"]\cr
[a[6] + a[5] + a[4] + a[3], "s.a11+a5+a3"]\cr
[a[6] + a[4] + 3 a[2] + 2 a[1], "s.a11+3a2"]\cr\footer\hss}}%
%[2 a[5] + a[4] + 2 a[2], "s.2a5+2a2+d5"]\cr\footer\hss}}%
\endinsert

\subsection{Classical Zariski pairs}\label{s.classical}
The list resulting from Table~\ref{tab.known}
contains a number of sextics of weight~$7$ with
at least two cusps. Perturbing a cusp of such a sextic,
%(inner or outer),
we obtain $30$ so called \emph{classical Zariski pairs},
\ie, pairs of irreducible sextics that share
the same set of singularities but differ by their Alexander
polynomials (see~\cite{JAG} for details and further references;
the sextic is/is not of torus type if the cusp
perturbed is, respectively, outer/inner). One can add to this list
the sets of singularities
$$
\bE_6\splus\bA_{11},\quad
\bE_6\splus\bA_8\splus\bA_2,\quad
\bA_{17},\quad
\bA_{11}\splus\bA_5,\quad
2\bA_8,
\eqtag\label{eq.CZP}
$$
which are realized by sextics with abelian fundamental groups in
Eyral, Oka~\cite{EyralOka.abelian}, thus obtaining $35$ classical
Zariski pairs. (Sextics of torus type realizing~\eqref{eq.CZP}
are constructed in this paper. The two other sets of
singularities discovered in~\cite{EyralOka.abelian} are already on
the list.) In each pair, the groups of the two curves are $\RBG3$
and $\CG6$.

According to~\cite{JAG} and~\cite{Aysegul}, each of the $35$ sets
of singularities obtained above is realized by exactly two
equisingular deformation families of irreducible sextics,
one of torus type and one not.
Altogether, there are $51$ classical Zariski pairs
of irreducible sextics (one of them
being, in fact, a triple). For all but one of them (the set of
singularities $(\bA_{14}\splus\bA_2)\splus2\bA_1$), the group of
the curve of torus type is known; it equals $\RBG3$.

\widestnumber\key{EO1}
\refstyle{C}

\enddocument